\documentclass[english,5p, final]{elsarticle}
\usepackage[T1]{fontenc}
\usepackage[latin9]{inputenc}
\usepackage{color}
\usepackage{amsmath}
\usepackage{amssymb}
\usepackage{esint}

\makeatletter
\renewcommand*{\appendixname}{}
\usepackage[colorlinks,breaklinks=true]{hyperref}
\usepackage[anythingbreaks]{breakurl}

\makeatother

\usepackage{babel}
\begin{document}

\begin{frontmatter}{}

\title{Exactly realizable desired trajectories}

\author{Jakob Löber}

\address{Institut für Theoretische Physik, Hardenbergstraße 36, EW 7-1, Technische
Universität Berlin, 10623 Berlin}
\begin{abstract}
Trajectory tracking of nonlinear dynamical systems with affine open-loop
controls is investigated. The control task is to enforce the system
state to follow a prescribed desired trajectory as closely as possible.
We introduce exactly realizable desired trajectories as these trajectories
which can be tracked exactly by an appropriate control. Exactly realizable
trajectories are characterized mathematically by means of Moore-Penrose
projectors constructed from the input matrix. The approach leads to
differential-algebraic systems of equations and is considerably simpler
than the related concept of system inversion. Furthermore, we identify
a particularly simple class of nonlinear affine control systems. Systems
in this class satisfy the so-called linearizing assumption and share
many properties with linear control systems. For example, conditions
for controllability can be formulated in terms of a rank condition
for a controllability matrix analogously to the Kalman rank condition
for linear time-invariant systems.
\end{abstract}

\end{frontmatter}{}

\section{Introduction}

A common approach to control is concerned with states as the object
to be controlled \cite{chen1995linear,khalil2002nonlinear}. Suppose
a controlled system, often called a plant in this context, has a certain
point $\boldsymbol{x}_{1}$ in state space, sometimes called the operating
point, at which the system works efficiently. The control task is
then to bring the system to the operating point $\boldsymbol{x}_{1}$,
and keep it there.\\
In contrast to that, here we develop an approach to control which
centers on the state trajectory over time, $\boldsymbol{x}\left(t\right)$,
as the object of interest. We distinguish between the controlled state
trajectory $\boldsymbol{x}\left(t\right)$ and the desired trajectory
$\boldsymbol{x}_{d}\left(t\right)$. The former is the trajectory
which the time-dependent state $\boldsymbol{x}\left(t\right)$ traces
out in state space under the action of a control signal, also called
an input signal. The latter is a fictitious reference trajectory for
the state over time. It is prescribed in analytical or numerical form
by the experimenter. Depending on the choice of the desired trajectory
$\boldsymbol{x}_{d}\left(t\right)$, the controlled state $\boldsymbol{x}\left(t\right)$
may or may not exactly follow $\boldsymbol{x}_{d}\left(t\right)$.\\
Of course, both approaches to control are closely related. A single
operating point in state space at which the system is to be kept is
nothing more than a degenerate state trajectory. Equivalently, any
state trajectory can be approximated by a succession of working points.\\
Trajectory tracking aims at enforcing, via a control signal $\boldsymbol{u}\left(t\right)$,
a system state $\boldsymbol{x}\left(t\right)$ to follow a prescribed
desired trajectory $\boldsymbol{x}_{d}\left(t\right)$ as closely
as possible within a time interval $t_{0}\leq t\leq t_{1}$. The distance
between $\boldsymbol{x}\left(t\right)$ and $\boldsymbol{x}_{d}\left(t\right)$
in function space can be measured by the functional 
\begin{align}
\mathcal{J} & =\frac{1}{2}\intop_{t_{0}}^{t_{1}}\text{d}t\left(\boldsymbol{x}\left(t\right)-\boldsymbol{x}_{d}\left(t\right)\right)^{2}.\label{eq:JDef}
\end{align}
The smallest possible value $\mathcal{J}=0$ is attained if and only
if the state $\boldsymbol{x}\left(t\right)$ follows the desired trajectory
exactly, i.e.,
\begin{align}
\boldsymbol{x}\left(t\right) & =\boldsymbol{x}_{d}\left(t\right)\label{eq:ExactlyRealizable}
\end{align}
for all times $t_{0}\leq t\leq t_{1}$. We call a desired trajectory
$\boldsymbol{x}_{d}\left(t\right)$ for which Eq. \eqref{eq:ExactlyRealizable}
holds an \textit{exactly realizable desired trajectory}. Clearly,
not every desired trajectory $\boldsymbol{x}_{d}\left(t\right)$ can
be exactly realized. The question addressed in this article is how,
for a given affine control system, exactly realizable desired trajectories
can be characterized mathematically.\textcolor{black}{}\\
\textcolor{black}{Tracking and regulation of desired outputs are common
problems in applications and have a long history of research. The
linear quadratic regulator \cite{bryson1969applied} is a cornerstone
of control theory. Further notable achievements are the solution of
the linear time-invariant (LTI) regulator problem by Francis \cite{francis1977linear}
and its generalization to nonlinear systems, the Byrnes-Isidori regulator
\cite{isidori1990output}. These regulators track desired trajectories
asymptotically and can deal with external disturbances and perturbations
of initial conditions. In contrast to that, here we consider the exact
tracking of desired trajectories in undisturbed systems by open-loop
control.}\\
\textcolor{black}{A concept closely related to} our work is that of
an inversion of control systems. The idea there is to find a second
controlled dynamical system which takes the desired output of the
original system as the input and outputs the input of the original
system. A control system is invertible when the corresponding input-output
map is injective. The resulting control signal is open-loop and often
referred to as feed-forward control. Stabilization of potential instabilities
can be accomplished by an additional feedback control. Early work
investigated the invertibility of LTI systems \cite{dorato1969inverse,silverman1969inversion,sain1969invertibility}.
Hirschorn analyzed invertibility of nonlinear systems for single \cite{hirschorn1979invertibilitya}
and multivariable \cite{hirschorn1979invertibility} input signals.
Inversion is commonly investigated by generating and analyzing a
hierarchy of auxiliary dynamical systems. Newer works focus on the
stability of inversion-based output tracking by combining system inversion
with feedback \cite{devasia1996nonlinear,zou2009optimal}.\\
The formalism necessary for the mathematical characterization of realizable
trajectories is introduced in Section \ref{sec:Formalism}. Section
\ref{sec:ExactlyRealizableTrajectories} defines the notion of exactly
realizable trajectories while Section \ref{sec:OutputTrajectoryRealizability}
introduces output realizability. Section \ref{sec:LinearizingAssumption}
proposes a basic assumption, called the linearizing assumption. This
assumption defines a class of nonlinear control systems which, to
a large extent, behave like linear systems. We demonstrate how controllability
can be recovered in our approach for systems satisfying the linearizing
assumption in Sections \ref{sec:Controllability} and \ref{sec:OutputControllability}.
Section \ref{sec:2Conclusions} concludes with a discussion and outlook.

\section{\label{sec:Formalism}Formalism}

Consider the affine control system for the state $\boldsymbol{x}\in\mathbb{R}^{n}$
with output $\boldsymbol{y}\in\mathbb{R}^{m}$ and control signal
$\boldsymbol{u}\in\mathbb{R}^{p}$,
\begin{align}
\boldsymbol{\dot{x}}\left(t\right) & =\boldsymbol{R}\left(\boldsymbol{x}\left(t\right)\right)+\boldsymbol{\mathcal{B}}\left(\boldsymbol{x}\left(t\right)\right)\boldsymbol{u}\left(t\right), & \boldsymbol{y}\left(t\right) & =\boldsymbol{h}\left(\boldsymbol{x}\left(t\right)\right).\label{eq:StateEq}
\end{align}
The time derivative is denoted by $\boldsymbol{\dot{x}}\left(t\right)=\frac{\text{d}}{\text{d}t}\boldsymbol{x}\left(t\right)$.
The dynamical system \eqref{eq:StateEq} is supplemented with the
initial condition $\boldsymbol{x}\left(t_{0}\right)=\boldsymbol{x}_{0}$.
The $n\times p$ input matrix $\boldsymbol{\mathcal{B}}\left(\boldsymbol{x}\right)$
may be state dependent and is assumed to have full rank,
\begin{align}
\text{rank}\left(\boldsymbol{\mathcal{B}}\left(\boldsymbol{x}\right)\right) & =p,
\end{align}
for all $\boldsymbol{x}\in\mathbb{R}^{n}$. The main elements of the
formalism introduced below are two complementary projection matrices
defined in terms of the input matrix $\boldsymbol{\mathcal{B}}\left(\boldsymbol{x}\right)$.
\newdefinition{defn}{Definition}
\begin{defn}The Moore-Penrose pseudo inverse \cite{CampbellJr.1991Generalized}
of $\boldsymbol{\mathcal{B}}\left(\boldsymbol{x}\right)$, denoted
by $\boldsymbol{\mathcal{B}}^{+}\left(\boldsymbol{x}\right)$, is
defined as the $p\times n$ matrix 
\begin{align}
\boldsymbol{\mathcal{B}}^{+}\left(\boldsymbol{x}\right) & =\left(\boldsymbol{\mathcal{B}}^{T}\left(\boldsymbol{x}\right)\boldsymbol{\mathcal{B}}\left(\boldsymbol{x}\right)\right)^{-1}\boldsymbol{\mathcal{B}}^{T}\left(\boldsymbol{x}\right).\label{eq:MoorePenrosePseudoInverse}
\end{align}
The Moore-Penrose projectors $\boldsymbol{\mathcal{P}}\left(\boldsymbol{x}\right)$
and $\boldsymbol{\mathcal{Q}}\left(\boldsymbol{x}\right)$ are $n\times n$
matrices defined as
\begin{align}
\boldsymbol{\mathcal{P}}\left(\boldsymbol{x}\right) & =\boldsymbol{\mathcal{B}}\left(\boldsymbol{x}\right)\boldsymbol{\mathcal{B}}^{+}\left(\boldsymbol{x}\right), & \boldsymbol{\mathcal{Q}}\left(\boldsymbol{x}\right) & =\boldsymbol{1}-\boldsymbol{\mathcal{P}}\left(\boldsymbol{x}\right).\label{eq:TimeDependentPQ}
\end{align}
\end{defn}
\newdefinition{rem}{Remark}
\begin{rem}Note that the $p\times p$ matrix $\boldsymbol{\mathcal{B}}^{T}\left(\boldsymbol{x}\right)\boldsymbol{\mathcal{B}}\left(\boldsymbol{x}\right)$
has full rank $p$ because $\boldsymbol{\mathcal{B}}\left(\boldsymbol{x}\right)$
is assumed to have full rank. Therefore, $\boldsymbol{\mathcal{B}}^{T}\left(\boldsymbol{x}\right)\boldsymbol{\mathcal{B}}\left(\boldsymbol{x}\right)$
is a quadratic symmetric non-singular matrix and its inverse exists.
From the definitions \eqref{eq:TimeDependentPQ} follow idempotence
\begin{align}
\boldsymbol{\mathcal{P}}\left(\boldsymbol{x}\right)\boldsymbol{\mathcal{P}}\left(\boldsymbol{x}\right) & =\boldsymbol{\mathcal{P}}\left(\boldsymbol{x}\right), & \boldsymbol{\mathcal{Q}}\left(\boldsymbol{x}\right)\boldsymbol{\mathcal{Q}}\left(\boldsymbol{x}\right) & =\boldsymbol{\mathcal{Q}}\left(\boldsymbol{x}\right),
\end{align}
and complementarity
\begin{align}
\boldsymbol{\mathcal{Q}}\left(\boldsymbol{x}\right)\boldsymbol{\mathcal{P}}\left(\boldsymbol{x}\right) & =\boldsymbol{\mathcal{P}}\left(\boldsymbol{x}\right)\boldsymbol{\mathcal{Q}}\left(\boldsymbol{x}\right)=\boldsymbol{0}.\label{eq:PTimesQ}
\end{align}
The projectors are symmetric,
\begin{align}
\boldsymbol{\mathcal{P}}^{T}\left(\boldsymbol{x}\right) & =\boldsymbol{\mathcal{P}}\left(\boldsymbol{x}\right), & \boldsymbol{\mathcal{Q}}^{T}\left(\boldsymbol{x}\right) & =\boldsymbol{\mathcal{Q}}\left(\boldsymbol{x}\right),
\end{align}
and their ranks are
\begin{align}
\text{rank}\left(\boldsymbol{\mathcal{P}}\left(\boldsymbol{x}\right)\right) & =p, & \text{rank}\left(\boldsymbol{\mathcal{Q}}\left(\boldsymbol{x}\right)\right) & =n-p,
\end{align}
independent of $\boldsymbol{x}$. Furthermore, multiplying $\boldsymbol{\mathcal{P}}\left(\boldsymbol{x}\right)$
and $\boldsymbol{\mathcal{Q}}\left(\boldsymbol{x}\right)$ from the
right with the input matrix $\boldsymbol{\mathcal{B}}\left(\boldsymbol{x}\right)$
yields the important relations 
\begin{align}
\boldsymbol{\mathcal{P}}\left(\boldsymbol{x}\right)\boldsymbol{\mathcal{B}}\left(\boldsymbol{x}\right) & =\boldsymbol{\mathcal{B}}\left(\boldsymbol{x}\right), & \boldsymbol{\mathcal{Q}}\left(\boldsymbol{x}\right)\boldsymbol{\mathcal{B}}\left(\boldsymbol{x}\right) & =\boldsymbol{0}.\label{eq:PTimesB}
\end{align}
Equation \eqref{eq:PTimesB} shows that the $p$ linearly independent
columns of $\boldsymbol{\mathcal{B}}\left(\boldsymbol{x}\right)$
are eigenvectors of $\boldsymbol{\mathcal{P}}\left(\boldsymbol{x}\right)$
to eigenvalue one and eigenvectors of $\boldsymbol{\mathcal{Q}}\left(\boldsymbol{x}\right)$
to eigenvalue zero. Alternatively, due to the idempotence of the projectors,
the $p$ eigenvectors of $\boldsymbol{\mathcal{P}}\left(\boldsymbol{x}\right)$
to eigenvalue one are given by the $j\in\left\{ 1,\dots,p\right\} $
linearly independent columns $\boldsymbol{p}_{j}\left(\boldsymbol{x}\right)$
of $\boldsymbol{\mathcal{P}}\left(\boldsymbol{x}\right)$. The remaining
$n-p$ eigenvectors are given by the $i\in\left\{ 1,\dots,n-p\right\} $
linearly independent columns $\boldsymbol{q}_{i}\left(\boldsymbol{x}\right)$
of $\boldsymbol{\mathcal{Q}}\left(\boldsymbol{x}\right)$.\end{rem}\begin{rem}Using
$\boldsymbol{1}=\boldsymbol{\mathcal{P}}\left(\boldsymbol{x}\right)+\boldsymbol{\mathcal{Q}}\left(\boldsymbol{x}\right)$,
any state vector $\boldsymbol{x}$ can be split up as
\begin{align}
\boldsymbol{x} & =\boldsymbol{\mathcal{P}}\left(\boldsymbol{x}\right)\boldsymbol{x}+\boldsymbol{\mathcal{Q}}\left(\boldsymbol{x}\right)\boldsymbol{x}=\boldsymbol{v}+\boldsymbol{w}.
\end{align}
Because $\boldsymbol{\mathcal{P}}\left(\boldsymbol{x}\right)$ has
rank $p$, only $p\leq n$ of the $n$ components of $\boldsymbol{v}=\boldsymbol{\mathcal{P}}\left(\boldsymbol{x}\right)\boldsymbol{x}$
are independent. Similarly, only $n-p$ components of $\boldsymbol{w}=\boldsymbol{\mathcal{Q}}\left(\boldsymbol{x}\right)\boldsymbol{x}$
are independent. $p$ independent components $\boldsymbol{\hat{v}}\in\mathbb{R}^{p}$
of $\boldsymbol{v}$ and $n-p$ independent components $\boldsymbol{\hat{w}}\in\mathbb{R}^{n-p}$
of $\boldsymbol{w}$ can be obtained as $\boldsymbol{\hat{v}}=\boldsymbol{\mathcal{\hat{P}}}^{T}\left(\boldsymbol{x}\right)\boldsymbol{x}$
and $\boldsymbol{\hat{w}}=\boldsymbol{\mathcal{\hat{Q}}}^{T}\left(\boldsymbol{x}\right)\boldsymbol{x}$,
respectively. Here, the $n\times p$ matrix $\boldsymbol{\mathcal{\hat{P}}}\left(\boldsymbol{x}\right)$
and the $n\times\left(n-p\right)$ matrix $\boldsymbol{\mathcal{\hat{Q}}}\left(\boldsymbol{x}\right)$
are constructed from the linearly independent columns $\boldsymbol{p}_{i}\left(\boldsymbol{x}\right)$
of $\boldsymbol{\mathcal{P}}\left(\boldsymbol{x}\right)$ and $\boldsymbol{q}_{i}\left(\boldsymbol{x}\right)$
of $\boldsymbol{\mathcal{Q}}\left(\boldsymbol{x}\right)$ as
\begin{align}
\boldsymbol{\mathcal{\hat{P}}}\left(\boldsymbol{x}\right) & =\left(\boldsymbol{p}_{1}\left(\boldsymbol{x}\right)|\dots|\boldsymbol{p}_{p}\left(\boldsymbol{x}\right)\right),\label{eq:PHat}\\
\boldsymbol{\mathcal{\hat{Q}}}\left(\boldsymbol{x}\right) & =\left(\boldsymbol{q}_{1}\left(\boldsymbol{x}\right)|\dots|\boldsymbol{q}_{n-p}\left(\boldsymbol{x}\right)\right).\label{eq:QHat}
\end{align}
If the projectors $\boldsymbol{\mathcal{P}}\left(\boldsymbol{x}\right)=\boldsymbol{\mathcal{P}}$
and $\boldsymbol{\mathcal{Q}}\left(\boldsymbol{x}\right)=\boldsymbol{\mathcal{Q}}$
are independent of the state $\boldsymbol{x}$, the vectors $\boldsymbol{v}=\boldsymbol{\mathcal{P}}\boldsymbol{x}$
and $\boldsymbol{w}=\boldsymbol{\mathcal{Q}}\boldsymbol{x}$ are simply
linear combinations of the original state components $\boldsymbol{x}$.
If $\boldsymbol{\mathcal{P}}\left(\boldsymbol{x}\right)$ and therefore
$\boldsymbol{\mathcal{Q}}\left(\boldsymbol{x}\right)=\boldsymbol{1}-\boldsymbol{\mathcal{P}}\left(\boldsymbol{x}\right)$
depends on $\boldsymbol{x}$, both vectors $\boldsymbol{v}$ and $\boldsymbol{w}$
are nonlinear functions of the state $\boldsymbol{x}$. However, a
state transformation can be found such that $\boldsymbol{v}$ and
$\boldsymbol{w}$ attain a particularly simple form. Being a projector,
$\boldsymbol{\mathcal{Q}}\left(\boldsymbol{x}\right)$ can be diagonalized
by a nonsingular $n\times n$ matrix $\boldsymbol{\mathcal{T}}\left(\boldsymbol{x}\right)$,
resulting in a diagonal $n\times n$ matrix $\boldsymbol{\mathcal{Q}}_{D}$,
\begin{align}
\boldsymbol{\mathcal{Q}}_{D} & =\boldsymbol{\mathcal{T}}{}^{-1}\left(\boldsymbol{x}\right)\boldsymbol{\mathcal{Q}}\left(\boldsymbol{x}\right)\boldsymbol{\mathcal{T}}\left(\boldsymbol{x}\right).
\end{align}
The first $p$ entries on the diagonal of $\boldsymbol{\mathcal{Q}}_{D}$
are zero while the remaining $n-p$ entries on the diagonal of $\boldsymbol{\mathcal{Q}}_{D}$
are one. The same matrix $\boldsymbol{\mathcal{T}}\left(\boldsymbol{x}\right)$
diagonalizes the projector $\boldsymbol{\mathcal{P}}\left(\boldsymbol{x}\right)$
as well. Defining the transformed state vector $\boldsymbol{\tilde{x}}$
as
\begin{align}
\boldsymbol{\tilde{x}} & =\boldsymbol{\mathcal{T}}^{-1}\left(\boldsymbol{x}\right)\boldsymbol{x},\label{eq:StateTransform}
\end{align}
the separation of the state can be seen to attain the particularly
simple form
\begin{align}
\boldsymbol{\tilde{x}} & =\boldsymbol{\mathcal{T}}^{-1}\left(\boldsymbol{x}\right)\boldsymbol{\mathcal{P}}\left(\boldsymbol{x}\right)\boldsymbol{x}+\boldsymbol{\mathcal{T}}^{-1}\left(\boldsymbol{x}\right)\boldsymbol{\mathcal{Q}}\left(\boldsymbol{x}\right)\boldsymbol{x}\nonumber \\
 & =\boldsymbol{\mathcal{P}}_{D}\boldsymbol{\tilde{x}}+\boldsymbol{\mathcal{Q}}_{D}\boldsymbol{\tilde{x}}=\boldsymbol{\tilde{v}}+\boldsymbol{\tilde{w}}
\end{align}
with
\begin{align}
\boldsymbol{\tilde{v}} & =\boldsymbol{\mathcal{P}}_{D}\boldsymbol{\tilde{x}}=\left(\tilde{x}_{1},\dots,\,\tilde{x}_{p},\,0,\dots,\,0\right)^{T},\\
\boldsymbol{\tilde{w}} & =\boldsymbol{\mathcal{Q}}_{D}\boldsymbol{\tilde{x}}=\left(0,\dots,\,0,\,\tilde{x}_{p+1},\dots,\,\tilde{x}_{n}\right)^{T}.
\end{align}
The matrix $\boldsymbol{\mathcal{T}}\left(\boldsymbol{x}\right)$
can be constructed from the eigenvectors of $\boldsymbol{\mathcal{Q}}\left(\boldsymbol{x}\right)$
in the usual manner. The state transformation $\boldsymbol{\mathcal{T}}^{-1}\left(\boldsymbol{x}\right)$
leads to a new affine control system for $\boldsymbol{\tilde{x}}$.
This new system may be viewed as a normal form of the affine control
system \eqref{eq:StateEq} \cite{loeber2016optimal}.\end{rem}In
the following theorem, we use the Moore-Penrose projectors to separate
the controlled state equation in two equations. The first one provides
a relation for the control signal, while the second equation is independent
of the control. \newtheorem{thm}{Theorem}
\begin{thm}\label{thm:Theorem0}Every affinely controlled state equation \eqref{eq:StateEq}
can be split in two separate equations. The first equation 
\begin{align}
\boldsymbol{\mathcal{Q}}\left(\boldsymbol{x}\right)\left(\boldsymbol{\dot{x}}-\boldsymbol{R}\left(\boldsymbol{x}\right)\right) & =\mathbf{0},\label{eq:ConstraintEq}
\end{align}
is independent of the control signal $\boldsymbol{u}$ and is called
the constraint equation. The second equation yields an expression
for the control $\boldsymbol{u}$ in terms of the controlled state
$\boldsymbol{x}$,

\begin{align}
\boldsymbol{u} & =\boldsymbol{\mathcal{B}}^{+}\left(\boldsymbol{x}\right)\left(\boldsymbol{\dot{x}}-\boldsymbol{R}\left(\boldsymbol{x}\right)\right).\label{eq:NonLinearSystemControlSolution}
\end{align}
\end{thm}
\newproof{pot1}{Proof of Theorem \ref{thm:Theorem0}}
\begin{pot1}The controlled state equation \eqref{eq:StateEq} can be written as
\begin{align}
\frac{d}{dt}\left(\boldsymbol{\mathcal{P}}\left(\boldsymbol{x}\right)\boldsymbol{x}+\boldsymbol{\mathcal{Q}}\left(\boldsymbol{x}\right)\boldsymbol{x}\right) & =\left(\boldsymbol{\mathcal{P}}\left(\boldsymbol{x}\right)+\boldsymbol{\mathcal{Q}}\left(\boldsymbol{x}\right)\right)\boldsymbol{R}\left(\boldsymbol{x}\right)\nonumber \\
 & +\left(\boldsymbol{\mathcal{P}}\left(\boldsymbol{x}\right)+\boldsymbol{\mathcal{Q}}\left(\boldsymbol{x}\right)\right)\boldsymbol{\mathcal{B}}\left(\boldsymbol{x}\right)\boldsymbol{u}.
\end{align}
Multiplying with $\boldsymbol{\mathcal{Q}}\left(\boldsymbol{x}\right)$
from the left and using Eq. \eqref{eq:PTimesB} yields Eq. \eqref{eq:ConstraintEq}.
Multiplying Eq. \eqref{eq:StateEq} by $\boldsymbol{\mathcal{B}}^{T}\left(\boldsymbol{x}\right)$
from the left yields
\begin{align}
\boldsymbol{\mathcal{B}}^{T}\left(\boldsymbol{x}\right)\boldsymbol{\dot{x}} & =\boldsymbol{\mathcal{B}}^{T}\left(\boldsymbol{x}\right)\boldsymbol{R}\left(\boldsymbol{x}\right)+\boldsymbol{\mathcal{B}}^{T}\left(\boldsymbol{x}\right)\boldsymbol{\mathcal{B}}\left(\boldsymbol{x}\right)\boldsymbol{u}.
\end{align}
Multiplying with $\left(\boldsymbol{\mathcal{B}}^{T}\left(\boldsymbol{x}\right)\boldsymbol{\mathcal{B}}\left(\boldsymbol{x}\right)\right)^{-1}$
from the left yields Eq. \eqref{eq:NonLinearSystemControlSolution}
for the control.\qed\end{pot1}

A formalism based on Moore-Penrose projectors can be introduced for
the output as well. We assume an output of the form
\begin{align}
\boldsymbol{y}\left(t\right) & =\boldsymbol{\mathcal{C}}\left(\boldsymbol{x}\left(t\right)\right)\boldsymbol{x}\left(t\right)+\boldsymbol{a}\left(t\right).
\end{align}
The vector $\boldsymbol{a}\left(t\right)$ is independent of the state
and can be absorbed in the function $\boldsymbol{y}\left(t\right)$.
We drop $\boldsymbol{a}\left(t\right)$ in the following. The $m\times n$
matrix $\boldsymbol{\mathcal{C}}\left(\boldsymbol{x}\right)$ with
$m\leq n$ is assumed to have full rank for all $\boldsymbol{x}$,
\begin{align}
\text{rank}\left(\boldsymbol{\mathcal{C}}\left(\boldsymbol{x}\right)\right) & =m.
\end{align}
\begin{defn}\label{def:Definition2}The Moore-Penrose pseudo inverse
of $\boldsymbol{\mathcal{C}}\left(\boldsymbol{x}\right)$, denoted
by $\boldsymbol{\mathcal{C}}^{+}\left(\boldsymbol{x}\right)$, is
the $n\times m$ matrix 
\begin{align}
\boldsymbol{\mathcal{C}}^{+}\left(\boldsymbol{x}\right) & =\boldsymbol{\mathcal{C}}^{T}\left(\boldsymbol{x}\right)\left(\boldsymbol{\mathcal{C}}\left(\boldsymbol{x}\right)\boldsymbol{\mathcal{C}}^{T}\left(\boldsymbol{x}\right)\right)^{-1}.
\end{align}
The Moore-Penrose projectors $\boldsymbol{\mathcal{M}}\left(\boldsymbol{x}\right)$
and $\boldsymbol{\mathcal{N}}\left(\boldsymbol{x}\right)$ are $n\times n$
matrices defined by
\begin{align}
\boldsymbol{\mathcal{M}}\left(\boldsymbol{x}\right) & =\boldsymbol{\mathcal{C}}^{+}\left(\boldsymbol{x}\right)\boldsymbol{\mathcal{C}}\left(\boldsymbol{x}\right), & \boldsymbol{\mathcal{N}}\left(\boldsymbol{x}\right) & =\boldsymbol{1}-\boldsymbol{\mathcal{M}}\left(\boldsymbol{x}\right).\label{eq:DefMN}
\end{align}
\end{defn}
\begin{rem}The ranks of $\boldsymbol{\mathcal{M}}\left(\boldsymbol{x}\right)$
and $\boldsymbol{\mathcal{N}}\left(\boldsymbol{x}\right)$ are
\begin{align}
\text{rank}\left(\boldsymbol{\mathcal{M}}\left(\boldsymbol{x}\right)\right) & =m, & \text{rank}\left(\boldsymbol{\mathcal{N}}\left(\boldsymbol{x}\right)\right) & =n-m,
\end{align}
and they satisfy 
\begin{align}
\boldsymbol{\mathcal{M}}\left(\boldsymbol{x}\right)\boldsymbol{\mathcal{C}}^{T}\left(\boldsymbol{x}\right) & =\boldsymbol{\mathcal{C}}^{T}\left(\boldsymbol{x}\right), & \boldsymbol{\mathcal{N}}\left(\boldsymbol{x}\right)\boldsymbol{\mathcal{C}}^{T}\left(\boldsymbol{x}\right) & =\boldsymbol{0}.
\end{align}
With the help of $\boldsymbol{\mathcal{M}}\left(\boldsymbol{x}\right)$
and $\boldsymbol{\mathcal{N}}\left(\boldsymbol{x}\right)$ the state
vector $\boldsymbol{x}\left(t\right)$ can be split up as 
\begin{alignat}{1}
\boldsymbol{x} & =\boldsymbol{\mathcal{M}}\left(\boldsymbol{x}\right)\boldsymbol{x}+\boldsymbol{\mathcal{N}}\left(\boldsymbol{x}\right)\boldsymbol{x}=\boldsymbol{\mathcal{C}}^{+}\left(\boldsymbol{x}\right)\boldsymbol{y}+\boldsymbol{\mathcal{N}}\left(\boldsymbol{x}\right)\boldsymbol{x}.\label{eq:SplitStateMN}
\end{alignat}
Thus the part $\boldsymbol{\mathcal{M}}\left(\boldsymbol{x}\right)\boldsymbol{x}$
can be expressed in terms of the output $\boldsymbol{y}$ while the
part $\boldsymbol{\mathcal{N}}\left(\boldsymbol{x}\right)\boldsymbol{x}$
remains undetermined.\end{rem}

\section{\label{sec:ExactlyRealizableTrajectories}Exactly realizable desired
trajectories}

Not every desired trajectory $\boldsymbol{x}_{d}$ can be realized
by control. Here, we formulate a condition which has to be satisfied
by a desired trajectory to be exactly realizable.
\begin{thm}\label{thm:Theorem1}The controlled state trajectory $\boldsymbol{x}\left(t\right)$
follows the desired trajectory $\boldsymbol{x}_{d}\left(t\right)$
exactly,
\begin{align}
\boldsymbol{x}\left(t\right) & =\boldsymbol{x}_{d}\left(t\right),
\end{align}
if and only if
\begin{enumerate}
\item \label{enu:RealizableTrajectory1}$\boldsymbol{x}_{d}\left(t\right)$
satisfies the constraint equation
\begin{align}
\boldsymbol{\mathcal{Q}}\left(\boldsymbol{x}_{d}\left(t\right)\right)\left(\boldsymbol{\dot{x}}_{d}\left(t\right)-\boldsymbol{R}\left(\boldsymbol{x}_{d}\left(t\right)\right)\right) & =\mathbf{0},\label{eq:ConstraintEquationForRealizableTrajectories}
\end{align}

\item \label{enu:RealizableTrajectory2}the initial value $\boldsymbol{x}_{d}\left(t_{0}\right)$
equals the initial value $\boldsymbol{x}_{0}$ of the controlled state
equation
\begin{align}
\boldsymbol{x}_{d}\left(t_{0}\right) & =\boldsymbol{x}_{0},\label{eq:InitialConditionForRealizableTrajectories}
\end{align}

\item the control signal enforcing $\boldsymbol{x}_{d}\left(t\right)$ is
given by
\begin{align}
\boldsymbol{u}\left(t\right) & =\boldsymbol{\mathcal{B}}^{+}\left(\boldsymbol{x}_{d}\left(t\right)\right)\left(\boldsymbol{\dot{x}}_{d}\left(t\right)-\boldsymbol{R}\left(\boldsymbol{x}_{d}\left(t\right)\right)\right).\label{eq:ControlSignal}
\end{align}

\end{enumerate}
\end{thm}\newproof{pot}{Proof of Theorem \ref{thm:Theorem1}}
\begin{pot}From Theorem \ref{thm:Theorem0} and $\boldsymbol{x}\left(t\right)=\boldsymbol{x}_{d}\left(t\right)$
follows the necessity of conditions \eqref{eq:ConstraintEquationForRealizableTrajectories},
\eqref{eq:InitialConditionForRealizableTrajectories} and \eqref{eq:ControlSignal}.
For sufficiency, expression \eqref{eq:ControlSignal} for the control
is used in the controlled state equation \eqref{eq:StateEq} to obtain
\begin{align}
\boldsymbol{\dot{x}} & =\boldsymbol{R}\left(\boldsymbol{x}\right)+\boldsymbol{\mathcal{B}}\left(\boldsymbol{x}\right)\boldsymbol{\mathcal{B}}^{+}\left(\boldsymbol{x}_{d}\right)\left(\boldsymbol{\dot{x}}_{d}-\boldsymbol{R}\left(\boldsymbol{x}_{d}\right)\right).\label{eq:E16}
\end{align}
Note that $\boldsymbol{\mathcal{B}}$ depends on the actual system
state $\boldsymbol{x}$ while $\boldsymbol{\mathcal{B}}^{+}$ depends
on the desired state $\boldsymbol{x}_{d}$. We introduce the difference
$\Delta\boldsymbol{x}\left(t\right)$ between true and desired state
as 
\begin{align}
\Delta\boldsymbol{x}\left(t\right) & =\boldsymbol{x}\left(t\right)-\boldsymbol{x}_{d}\left(t\right).
\end{align}
Using Eq. \eqref{eq:E16}, the ODE for $\Delta\boldsymbol{x}$ reads
\begin{align}
\Delta\boldsymbol{\dot{x}} & =\boldsymbol{R}\left(\Delta\boldsymbol{x}+\boldsymbol{x}_{d}\right)-\boldsymbol{\dot{x}}_{d}\nonumber \\
 & +\boldsymbol{\mathcal{B}}\left(\Delta\boldsymbol{x}+\boldsymbol{x}_{d}\right)\boldsymbol{\mathcal{B}}^{+}\left(\boldsymbol{x}_{d}\right)\left(\boldsymbol{\dot{x}}_{d}-\boldsymbol{R}\left(\boldsymbol{x}_{d}\right)\right),\label{eq:ControlledStateEquationForDeltax}\\
\Delta\boldsymbol{x}\left(t_{0}\right) & =\boldsymbol{x}\left(t_{0}\right)-\boldsymbol{x}_{d}\left(t_{0}\right).
\end{align}
If the desired trajectory $\boldsymbol{x}_{d}$ satisfies initially
Eq. \eqref{eq:InitialConditionForRealizableTrajectories}, the initial
condition for Eq. \eqref{eq:ControlledStateEquationForDeltax} vanishes,
$\Delta\boldsymbol{x}\left(t_{0}\right)=\boldsymbol{0}$. If condition
\eqref{eq:ConstraintEquationForRealizableTrajectories} is satisfied,
then $\Delta\boldsymbol{x}\left(t\right)=\mathbf{0}$ is a stationary
point of Eq. \eqref{eq:ControlledStateEquationForDeltax}, 
\begin{align}
\Delta\boldsymbol{\dot{x}} & =\boldsymbol{R}\left(\boldsymbol{x}_{d}\right)-\boldsymbol{\dot{x}}_{d}+\boldsymbol{\mathcal{P}}\left(\boldsymbol{x}_{d}\right)\left(\boldsymbol{\dot{x}}_{d}-\boldsymbol{R}\left(\boldsymbol{x}_{d}\right)\right)\nonumber \\
 & =\boldsymbol{\mathcal{Q}}\left(\boldsymbol{x}_{d}\right)\left(\boldsymbol{\dot{x}}_{d}-\boldsymbol{R}\left(\boldsymbol{x}_{d}\right)\right)=\boldsymbol{0},
\end{align}
and so $\boldsymbol{x}\left(t\right)=\boldsymbol{x}_{d}\left(t\right)$
remains a solution to Eq. \eqref{eq:E16} for all times.\qed\end{pot}
\begin{rem}Because the control signal $\boldsymbol{u}\left(t\right)$ consists
of $p\leq n$ independent components, at most $p$ one-to-one relations
between state components and control components can be found. Thus,
maximally $p$ components of $\boldsymbol{x}_{d}\left(t\right)$ can
be prescribed by the experimenter, while the remaining $n-p$ components
are free. The time evolution of these $n-p$ components is fixed by
the constraint equation \eqref{eq:ConstraintEquationForRealizableTrajectories}.
This motivates the name constraint equation: for an arbitrary desired
trajectory $\boldsymbol{x}_{d}\left(t\right)$ to be exactly realizable,
it is constrained by Eq. \eqref{eq:ConstraintEquationForRealizableTrajectories}.
One possibility to obtain $n-p$ linearly independent equations from
the constraint equation is to use the matrix $\boldsymbol{\mathcal{\hat{Q}}}$
defined in Eq. \eqref{eq:QHat} as 
\begin{align}
\boldsymbol{\mathcal{\hat{Q}}}\left(\boldsymbol{x}_{d}\right)\left(\boldsymbol{\dot{x}}_{d}-\boldsymbol{R}\left(\boldsymbol{x}_{d}\right)\right) & =\mathbf{0}.
\end{align}

\end{rem}
\begin{rem}The necessity to satisfy the initial conditions $\boldsymbol{x}_{d}\left(t_{0}\right)=\boldsymbol{x}\left(t_{0}\right)=\boldsymbol{x}_{0}$,
Eq. \eqref{eq:InitialConditionForRealizableTrajectories}, leaves
us with two possibilities. Either the initial state $\boldsymbol{x}\left(t_{0}\right)=\boldsymbol{x}_{0}$
of the system can be prepared such that it equals the initial value
$\boldsymbol{x}_{d}\left(t_{0}\right)$ of a given desired trajectory
$\boldsymbol{x}_{d}$. Or the desired trajectory $\boldsymbol{x}_{d}$
must be designed such that it starts from the observed initial state
$\boldsymbol{x}_{0}$ of the system. The control signal as given by
Eq. \eqref{eq:ControlSignal} does neither depend on the current nor
on the previous state $\boldsymbol{x}$ of the system and is an open-loop
control signal. Only the initial state $\boldsymbol{x}_{0}$ of the
actual controlled system enters via the initial condition Eq. \eqref{eq:InitialConditionForRealizableTrajectories}
for the constraint equation. In general, the controlled system may
suffer from instability. For example, it is impossible to prepare
a real systems exactly in the initial state $\boldsymbol{x}_{0}$.
Furthermore, a mathematical model must be viewed as an approximation
of a real system. Disturbances which are not taken into account in
the model affect the time evolution of the state. An additional feedback
control may be sufficient to stabilize unstable desired trajectories.
However, a thorough discussion of these issues is outside the scope
of this article. \end{rem}

Equation \eqref{eq:ControlSignal} employs the Moore-Penrose pseudo
inverse $\boldsymbol{\mathcal{B}}^{+}\left(\boldsymbol{x}\right)$
to obtain the control signal in terms of the desired trajectory. This
choice for a generalized inverse matrix is not unique. Note that any
$p\times n$ matrix $\boldsymbol{\mathcal{B}}^{g}$ satisfying $\boldsymbol{\mathcal{B}}\boldsymbol{\mathcal{B}}^{g}\boldsymbol{\mathcal{B}}=\boldsymbol{\mathcal{B}}$
is called a generalized inverse of the $n\times p$ matrix $\boldsymbol{\mathcal{B}}$.
Indeed, any $p\times n$ matrix $\boldsymbol{\mathcal{K}}\left(\boldsymbol{x}\right)$
with the property
\begin{align}
\text{rank}\left(\boldsymbol{\mathcal{K}}\left(\boldsymbol{x}\right)\boldsymbol{\mathcal{B}}\left(\boldsymbol{x}\right)\right) & =p
\end{align}
for all $\boldsymbol{x}$ can be used to construct a generalized inverse
$p\times n$ matrix $\boldsymbol{\mathcal{B}}^{g}\left(\boldsymbol{x}\right)$
of $\boldsymbol{\mathcal{B}}\left(\boldsymbol{x}\right)$ as
\begin{align}
\boldsymbol{\mathcal{B}}^{g}\left(\boldsymbol{x}\right) & =\left(\boldsymbol{\mathcal{K}}\left(\boldsymbol{x}\right)\boldsymbol{\mathcal{B}}\left(\boldsymbol{x}\right)\right)^{-1}\boldsymbol{\mathcal{K}}\left(\boldsymbol{x}\right).\label{eq:GeneralizedInverse}
\end{align}
With the choice $\boldsymbol{\mathcal{K}}\left(\boldsymbol{x}\right)=\boldsymbol{\mathcal{B}}^{T}\left(\boldsymbol{x}\right)$,
$\boldsymbol{\mathcal{B}}^{g}\left(\boldsymbol{x}\right)$ becomes
the Moore-Penrose pseudo inverse $\boldsymbol{\mathcal{B}}^{+}\left(\boldsymbol{x}\right)$.
We demonstrate the uniqueness of the control solution Eq. \eqref{eq:ControlSignal}
and its independence of the choice of $\boldsymbol{\mathcal{K}}\left(\boldsymbol{x}\right)$
in the following theorem. In principle, any generalized inverse constructed
as in Eq. \eqref{eq:GeneralizedInverse} may be used to formulate
the constraint equation and the control solution.
\begin{thm}\label{thm:Theorem1.5}Let $\boldsymbol{\mathcal{K}}_{i}\left(\boldsymbol{x}\right)$
with $i\in\left\{ 1,2\right\} $ be two $p\times n$ matrices with
the property 
\begin{align}
\text{\textup{rank}}\left(\boldsymbol{\mathcal{K}}_{i}\left(\boldsymbol{x}\right)\boldsymbol{\mathcal{B}}\left(\boldsymbol{x}\right)\right) & =p.
\end{align}
Define two generalized inverses as
\begin{align}
\boldsymbol{\mathcal{B}}_{i}^{g}\left(\boldsymbol{x}\right) & =\left(\boldsymbol{\mathcal{K}}_{i}\left(\boldsymbol{x}\right)\boldsymbol{\mathcal{B}}\left(\boldsymbol{x}\right)\right)^{-1}\boldsymbol{\mathcal{K}}_{i}\left(\boldsymbol{x}\right),
\end{align}
and corresponding projectors
\begin{align}
\boldsymbol{\mathcal{P}}_{i}^{g}\left(\boldsymbol{x}\right) & =\boldsymbol{\mathcal{B}}\left(\boldsymbol{x}\right)\boldsymbol{\mathcal{B}}_{i}^{g}\left(\boldsymbol{x}\right), & \boldsymbol{\mathcal{Q}}_{i}^{g}\left(\boldsymbol{x}\right) & =\boldsymbol{1}-\boldsymbol{\mathcal{P}}_{i}^{g}\left(\boldsymbol{x}\right).
\end{align}
The control signals expressed in terms of the desired trajectory are
\begin{align}
\boldsymbol{u}_{i}\left(t\right) & =\boldsymbol{\mathcal{B}}_{i}^{g}\left(\boldsymbol{x}_{d}\left(t\right)\right)\left(\boldsymbol{\dot{x}}_{d}\left(t\right)-\boldsymbol{R}\left(\boldsymbol{x}_{d}\left(t\right)\right)\right),
\end{align}
with desired trajectory $\boldsymbol{x}_{d}$ constrained by
\begin{align}
\boldsymbol{0} & =\boldsymbol{\mathcal{Q}}_{i}^{g}\left(\boldsymbol{x}_{d}\left(t\right)\right)\left(\boldsymbol{\dot{x}}_{d}\left(t\right)-\boldsymbol{R}\left(\boldsymbol{x}_{d}\left(t\right)\right)\right).
\end{align}
Then both control signals are identical, 
\begin{align}
\boldsymbol{u}_{1}\left(t\right) & =\boldsymbol{u}_{2}\left(t\right).
\end{align}

\end{thm}\newproof{pot1.5}{Proof of Theorem \ref{thm:Theorem1.5}}
\begin{pot1.5}Multiplying the difference $\boldsymbol{u}_{1}-\boldsymbol{u}_{2}$
by $\boldsymbol{\mathcal{B}}\left(\boldsymbol{x}_{d}\right)$ and
exploiting the definitions of the projectors as well as the constraint
equations yields 
\begin{align}
 & \boldsymbol{\mathcal{B}}\left(\boldsymbol{x}_{d}\right)\left(\boldsymbol{u}_{1}-\boldsymbol{u}_{2}\right)\nonumber \\
= & \left(\boldsymbol{\mathcal{B}}\left(\boldsymbol{x}_{d}\right)\boldsymbol{\mathcal{B}}_{1}^{g}\left(\boldsymbol{x}_{d}\right)-\boldsymbol{\mathcal{B}}\left(\boldsymbol{x}_{d}\right)\boldsymbol{\mathcal{B}}_{2}^{g}\left(\boldsymbol{x}_{d}\right)\right)\left(\boldsymbol{\dot{x}}_{d}-\boldsymbol{R}\left(\boldsymbol{x}_{d}\right)\right)\nonumber \\
= & \left(\boldsymbol{\mathcal{P}}_{1}^{g}\left(\boldsymbol{x}_{d}\right)-\boldsymbol{\mathcal{P}}_{2}^{g}\left(\boldsymbol{x}_{d}\right)\right)\left(\boldsymbol{\dot{x}}_{d}-\boldsymbol{R}\left(\boldsymbol{x}_{d}\right)\right)\nonumber \\
= & \left(\boldsymbol{\mathcal{Q}}_{2}^{g}\left(\boldsymbol{x}_{d}\right)-\boldsymbol{\mathcal{Q}}_{1}^{g}\left(\boldsymbol{x}_{d}\right)\right)\left(\boldsymbol{\dot{x}}_{d}-\boldsymbol{R}\left(\boldsymbol{x}_{d}\right)\right)=\boldsymbol{0}.
\end{align}
Thus $\boldsymbol{u}_{1}-\boldsymbol{u}_{2}$ lies in the null space
of the input matrix $\boldsymbol{\mathcal{B}}\left(\boldsymbol{x}_{d}\right)$.
Because $\boldsymbol{\mathcal{B}}\left(\boldsymbol{x}_{d}\right)$
has full rank $p$ by assumption, its null space contains only $\boldsymbol{0}$
and so $\boldsymbol{u}_{1}\left(t\right)=\boldsymbol{u}_{2}\left(t\right)$.\qed
\end{pot1.5}

\section{\label{sec:OutputTrajectoryRealizability}Output realizability}

The constraint equation does not dictate which state components are
prescribed and which are fixed by the constraint equation. In general,
we have the freedom to choose an output $\boldsymbol{y}=\boldsymbol{h}\left(\boldsymbol{x}\right)$
with $m$ components such that under the action of control, a prescribed
desired output $\boldsymbol{y}_{d}\left(t\right)$ is exactly realized,
$\boldsymbol{y}_{d}\left(t\right)=\boldsymbol{y}\left(t\right)$.
\begin{thm}\label{thm:Theorem1.6}Let $\boldsymbol{x}_{d}\left(t\right)$ be
an exactly realizable trajectory, i.e., it satisfies the constraint
equation \eqref{eq:ConstraintEquationForRealizableTrajectories} and
the initial condition Eq. \eqref{eq:InitialConditionForRealizableTrajectories}.
If $\boldsymbol{x}_{d}\left(t\right)$ additionally satisfies 
\begin{align}
\boldsymbol{y}_{d}\left(t\right) & =\boldsymbol{h}\left(\boldsymbol{x}_{d}\left(t\right)\right),\label{eq:OutputRelation}
\end{align}
then the output $\boldsymbol{y}_{d}\left(t\right)$ is realized exactly,
i.e.,
\begin{align}
\boldsymbol{y}\left(t\right) & =\boldsymbol{y}_{d}\left(t\right).
\end{align}

\end{thm}\newproof{pot1.6}{Proof of Theorem \ref{thm:Theorem1.6}}
\begin{pot1.6}Let $\Delta\boldsymbol{y}\left(t\right)$ be defined as
\begin{align}
\Delta\boldsymbol{y}\left(t\right) & =\boldsymbol{y}_{d}\left(t\right)-\boldsymbol{y}\left(t\right)=\boldsymbol{h}\left(\boldsymbol{x}_{d}\left(t\right)\right)-\boldsymbol{h}\left(\boldsymbol{x}\left(t\right)\right)\nonumber \\
 & =\boldsymbol{h}\left(\boldsymbol{x}_{d}\left(t\right)\right)-\boldsymbol{h}\left(\Delta\boldsymbol{x}\left(t\right)+\boldsymbol{x}_{d}\left(t\right)\right),
\end{align}
with $\Delta\boldsymbol{x}\left(t\right)=\boldsymbol{x}\left(t\right)-\boldsymbol{x}_{d}\left(t\right)$.
For an exactly realizable trajectory we have $\Delta\boldsymbol{x}\left(t\right)=\boldsymbol{0}$
and so $\boldsymbol{y}\left(t\right)=\boldsymbol{y}_{d}\left(t\right)$.\qed
\end{pot1.6}

The solution to a control problem consists of a solution $\boldsymbol{x}$
to the controlled state equation \eqref{eq:StateEq} and a solution
$\boldsymbol{u}$ for the control signal. Within the framework of
exactly realizable trajectories, these are given by
\begin{align}
\boldsymbol{u} & =\boldsymbol{\mathcal{B}}^{+}\left(\boldsymbol{x}_{d}\right)\left(\boldsymbol{\dot{x}}_{d}-\boldsymbol{R}\left(\boldsymbol{x}_{d}\right)\right), & \boldsymbol{x} & =\boldsymbol{x}_{d}.
\end{align}
For $\boldsymbol{x}=\boldsymbol{x}_{d}$ to hold, $\boldsymbol{x}_{d}$
must satisfy the constraint equation \eqref{eq:ConstraintEquationForRealizableTrajectories}.
For output realizability, $\boldsymbol{x}_{d}$ additionally has to
satisfy the output relation \eqref{eq:OutputRelation}. Thus the only
equations which remain to be solved is the system of $n-p+m$ inhomogeneous
differential-algebraic equations (DAE) given by 
\begin{align}
\boldsymbol{y}_{d} & =\boldsymbol{h}\left(\boldsymbol{x}_{d}\right), & \boldsymbol{\mathcal{Q}}\left(\boldsymbol{x}_{d}\right)\left(\boldsymbol{\dot{x}}_{d}-\boldsymbol{R}\left(\boldsymbol{x}_{d}\right)\right) & =\mathbf{0}.\label{eq:OutputConstraintDAE}
\end{align}
The desired output $\boldsymbol{y}_{d}\left(t\right)$ represents
an inhomogeneity and renders Eqs. \eqref{eq:OutputConstraintDAE}
a non-autonomous DAE. Equations \eqref{eq:OutputConstraintDAE} have
to be solved for $n-p+m$ components of $\boldsymbol{x}_{d}$ together
with the $n$ initial conditions $\boldsymbol{x}_{d}\left(t_{0}\right)=\boldsymbol{x}_{0}$.
In principle, solutions can exist as long as $m\leq p$. If $m<p$,
$p-m$ components of $\boldsymbol{x}_{d}$ may be freely chosen as
long as they do not violate the initial conditions. We expect a solution
$\boldsymbol{x}_{d}\left(t\right)$ to Eqs. \eqref{eq:OutputConstraintDAE}
to depend on the entire history of the output $\boldsymbol{y}_{d}$
from the initial time $t_{0}$ up to the current time $t$. In a final
step, the solution for $\boldsymbol{x}_{d}$ obtained from Eqs. \eqref{eq:OutputConstraintDAE}
is used in the control signal $\boldsymbol{u}$, Eq. \eqref{eq:ControlSignal},
to eliminate $n-p+m$ components of $\boldsymbol{x}_{d}$. In general,
$\boldsymbol{u}\left(t\right)$ depends on the entire history of $\boldsymbol{y}_{d}$
up to the current time $t$.

Being a system of DAEs, Eqs. \eqref{eq:OutputConstraintDAE} cannot
accommodate all $n$ initial conditions. Consequences of the initial
conditions can be distinguished as follows. First, evaluating the
output relation at $t=t_{0}$ imposes $m$ conditions on the desired
output as
\begin{align}
\boldsymbol{y}_{d}\left(t_{0}\right) & =\boldsymbol{h}\left(\boldsymbol{x}_{0}\right).
\end{align}
Second, $r$ initial conditions are accommodated by the constants
of integration arising in the constraint equation. Third, in case
that $n-m-r=l>0$, $l$ additional relations between $\boldsymbol{y}_{d}\left(t_{0}\right)$
and $\boldsymbol{x}_{0}$ must be satisfied. The latter conditions
also involve the time derivative of $\boldsymbol{y}_{d}$ at the initial
time $t=t_{0}$. They ensure that the constraint equation is satisfied
also at $t=t_{0}$,
\begin{align}
\boldsymbol{\mathcal{Q}}\left(\boldsymbol{x}_{d}\left(t_{0}\right)\right)\left(\boldsymbol{\dot{x}}_{d}\left(t_{0}\right)-\boldsymbol{R}\left(\boldsymbol{x}_{d}\left(t_{0}\right)\right)\right) & =\mathbf{0}.
\end{align}
We illustrate output realizability with the help of two examples.
\newdefinition{example}{Example}
\begin{example}\label{exp:FHN}Consider the system
\begin{align}
\dot{x}_{1}\left(t\right) & =x_{2}\left(t\right),\\
\dot{x}_{2}\left(t\right) & =R\left(x_{1}\left(t\right),x_{2}\left(t\right)\right)+B\left(x_{1}\left(t\right),x_{2}\left(t\right)\right)u\left(t\right),
\end{align}
with state vector $\boldsymbol{x}=\left(x_{1},x_{2}\right)^{T}$,
nonlinearity $\boldsymbol{R}\left(\boldsymbol{x}\right)=\left(x_{2},R\left(x_{1},x_{2}\right)\right)^{T}$,
and input matrix $\boldsymbol{\mathcal{B}}\left(\boldsymbol{x}\right)=\left(0,B\left(x_{1},x_{2}\right)\right)^{T}$.
The initial condition is $\boldsymbol{x}\left(t_{0}\right)=\boldsymbol{x}_{0}=\left(x_{1,0},x_{2,0}\right)^{T}$.
The assumption of full rank of $\boldsymbol{\mathcal{B}}\left(\boldsymbol{x}\right)$
for all $\boldsymbol{x}$ implies $B\left(x_{1},x_{2}\right)\neq0$.
The system represents Newton's equation of motion for a point particle
with unit mass. The particle moves with position $x_{1}$ and velocity
$x_{2}$ in one spatial dimensional under the influence of an external
force $R$ and a control force $Bu$. The Moore-Penrose pseudo inverse
of the input matrix $\boldsymbol{\mathcal{B}}$ and the corresponding
projectors are
\begin{align}
\boldsymbol{\mathcal{B}}^{+}\left(\boldsymbol{x}\right) & =B\left(x_{1},x_{2}\right)^{-2}\left(0,B\left(x_{1},x_{2}\right)\right),\\
\boldsymbol{\mathcal{P}}\left(\boldsymbol{x}\right) & =\boldsymbol{\mathcal{P}}=\left(\begin{array}{cc}
0 & 0\\
0 & 1
\end{array}\right),\\
\boldsymbol{\mathcal{Q}}\left(\boldsymbol{x}\right) & =\boldsymbol{\mathcal{Q}}=\mathbf{1}-\boldsymbol{\mathcal{P}}=\left(\begin{array}{cc}
1 & 0\\
0 & 0
\end{array}\right).
\end{align}
The constraint equation and the control signal are
\begin{align}
\dot{x}_{1,d}\left(t\right) & =x_{2,d}\left(t\right), & u\left(t\right) & =\frac{\dot{x}_{2,d}\left(t\right)-R\left(x_{1,d}\left(t\right),x_{2,d}\left(t\right)\right)}{B\left(x_{1,d}\left(t\right),x_{2,d}\left(t\right)\right)},
\end{align}
respectively.

Let us first assume that the output $y$ is given by the velocity
$x_{2}$, $y\left(t\right)=x_{2}\left(t\right)$. The control task
is then to enforce a velocity over time prescribed by the experimenter
in form of the desired output $y_{d}$. The system of DAEs for output
realizability \eqref{eq:OutputConstraintDAE} becomes
\begin{align}
x_{2,d}\left(t\right) & =y_{d}\left(t\right), & \dot{x}_{1,d}\left(t\right) & =x_{2,d}\left(t\right).
\end{align}
The constraint equation is a differential equation for the position
$x_{1,d}$ and yields
\begin{align}
x_{1,d}\left(t\right) & =x_{1,d}\left(t_{0}\right)+\int_{t_{0}}^{t}\text{d}\tau y_{d}\left(\tau\right).
\end{align}
To satisfy the initial condition $\boldsymbol{x}_{d}\left(t_{0}\right)=\boldsymbol{x}_{0}$,
we must have
\begin{align}
x_{1,d}\left(t_{0}\right) & =x_{1,0}, & x_{2,d}\left(t_{0}\right) & =x_{2,0}=y_{d}\left(t_{0}\right).
\end{align}
Thus, for a given desired output $y_{d}$, the system has to be prepared
in the state $\boldsymbol{x}_{0}=\left(x_{1,0},y_{d}\left(t_{0}\right)\right)^{T}$,
with $x_{1,0}$ a free parameter. Or, for a given system with initial
state $\boldsymbol{x}_{0}$, the desired output must be chosen such
that initially $y_{d}\left(t_{0}\right)=x_{2,0}$. The control signal
becomes an expression depending only on the desired output $y_{d}$
and the initial value $x_{1,0}$,
\begin{align}
u\left(t\right) & =\frac{\dot{y}_{d}\left(t\right)-R\left(x_{1,0}+\intop_{t_{0}}^{t}\text{d}\tau y_{d}\left(\tau\right),y_{d}\left(t\right)\right)}{B\left(x_{1,0}+\intop_{t_{0}}^{t}\text{d}\tau y_{d}\left(\tau\right),y_{d}\left(t\right)\right)}.
\end{align}
The context of a mechanical control system allows the following interpretation
of our approach. The constraint equation can be viewed as the definition
of the velocity of a point particle. Neither an external force $R$
nor a control force $Bu$ can change that definition. With only a
single control force, position $x_{1}$ and velocity $x_{2}$ over
time cannot be controlled independently from each other.

Instead of choosing the velocity $x_{2}$ as the output, we may choose
the position $y\left(t\right)=x_{1}\left(t\right)$ as well. The DAE
\eqref{eq:OutputConstraintDAE} becomes
\begin{align}
x_{1,d}\left(t\right) & =y_{d}\left(t\right), & \dot{x}_{1,d}\left(t\right) & =x_{2,d}\left(t\right).
\end{align}
Here, the constraint equation is an algebraic equation for the desired
velocity $x_{2,d}$, $x_{2,d}\left(t\right)=\dot{y}_{d}\left(t\right)$,
which is used eliminate $x_{2,d}\left(t\right)$ from the control
signal. We obtain an expression which depends on the desired output
only,
\begin{align}
u\left(t\right) & =\frac{\ddot{y}_{d}\left(t\right)-R\left(y_{d}\left(t\right),\dot{y}_{d}\left(t\right)\right)}{B\left(y_{d}\left(t\right),\dot{y}_{d}\left(t\right)\right)}.\label{eq:ComputedTorque}
\end{align}
To satisfy the initial condition $\boldsymbol{x}_{d}\left(t_{0}\right)=\boldsymbol{x}_{0}$,
the desired output $y_{d}\left(t\right)$ must satisfy the two initial
conditions 
\begin{align}
x_{1,d}\left(t_{0}\right) & =x_{1,0}=y_{d}\left(t_{0}\right), & x_{2,d}\left(t_{0}\right) & =x_{2,0}=\dot{y}_{d}\left(t_{0}\right).
\end{align}
In this case, the initial desired output has to satisfy two conditions
because the constraint equation is a purely algebraic equation which
does not allow for an initial condition. Equation \eqref{eq:ComputedTorque}
is known as a variant of the so-called computed torque formula and
has found widespread application in robotics \cite{lewis1993control,deWit1996theory,angeles2002fundamentals}.
\end{example}
\begin{example}We consider the example from \cite{devasia1996nonlinear} with $n=4$
and $p=1$. The output is $y=x_{1}-3x_{3}$, and the nonlinearity
$\boldsymbol{R}$ and input matrix $\boldsymbol{\mathcal{B}}$ are
given by 
\begin{align}
\boldsymbol{R}\left(\boldsymbol{x}\right) & =\left(\begin{array}{c}
x_{2}-x_{1}\\
x_{1}^{3}-3x_{2}\\
x_{1}-2x_{3}\\
x_{3}^{2}-x_{4}
\end{array}\right), & \boldsymbol{\mathcal{B}}\left(\boldsymbol{x}\right) & =\left(\begin{array}{c}
0\\
2+\sin^{2}\left(x_{4}\right)\\
0\\
0
\end{array}\right).
\end{align}
The control signal in terms of the desired trajectory is
\begin{align}
u\left(t\right) & =\frac{\dot{x}_{2,d}\left(t\right)-x_{1,d}^{3}\left(t\right)+3x_{2,d}\left(t\right)}{2+\sin^{2}\left(x_{4,d}\left(t\right)\right)}.
\end{align}
Equations \eqref{eq:OutputConstraintDAE} becomes 
\begin{align}
y_{d}\left(t\right) & =x_{1,d}\left(t\right)-3x_{3,d}\left(t\right),\\
\dot{x}_{1,d}\left(t\right) & =-x_{1,d}\left(t\right)+x_{2,d}\left(t\right),\\
\dot{x}_{3,d}\left(t\right) & =x_{1,d}\left(t\right)-2x_{3,d}\left(t\right),\\
\dot{x}_{4,d}\left(t\right) & =x_{3,d}^{2}\left(t\right)-x_{4,d}\left(t\right),
\end{align}
which must be solved for $\boldsymbol{x}_{d}\left(t\right)$. We obtain
\begin{align}
x_{1,d}\left(t\right) & =3\int_{t_{0}}^{t}\text{d}\tau\text{e}^{t-\tau}y_{d}\left(\tau\right)+y_{d}\left(t\right)\nonumber \\
 & +\text{e}^{t-t_{0}}\left(x_{1,d}\left(t_{0}\right)-y_{d}\left(t_{0}\right)\right),\\
x_{2,d}\left(t\right) & =\dot{x}_{1,d}\left(t\right)+x_{1,d}\left(t\right),\\
x_{3,d}\left(t\right) & =\dfrac{1}{3}\left(x_{1,d}\left(t\right)-y_{d}\left(t\right)\right),\\
x_{4,d}\left(t\right) & =\text{e}^{t_{0}-t}x_{4}\left(t_{0}\right)+\frac{1}{9}\int_{t_{0}}^{t}\text{d}\tau\text{e}^{\tau+t-2t_{0}}\left(f\left(\tau\right)\right)^{2},
\end{align}
with
\begin{align}
f\left(t\right) & =x_{1}\left(t_{0}\right)-y_{d}\left(t_{0}\right)+3\int_{t_{0}}^{t}\text{d}\tau\text{e}^{t_{0}-\tau}y_{d}\left(\tau\right).
\end{align}
Enforcing the initial condition $\boldsymbol{x}_{d}\left(t_{0}\right)=\boldsymbol{x}_{0}$
yields two additional relations for the initial desired output,
\begin{align}
x_{2,d}\left(t_{0}\right) & =x_{2,0}=2\left(x_{1,0}+y_{d}\left(t_{0}\right)\right)+\dot{y}_{d}\left(t_{0}\right),\\
x_{3,d}\left(t_{0}\right) & =x_{3,0}=\frac{1}{3}\left(x_{1,0}-y_{d}\left(t_{0}\right)\right),
\end{align}
while two values of the initial conditions are free parameters,
\begin{align}
x_{1,d}\left(t_{0}\right) & =x_{1,0}, & x_{4,d}\left(t_{0}\right) & =x_{4,0}.
\end{align}

\end{example}

\section{\label{sec:LinearizingAssumption}Linearizing assumption}

Instead of having to solve the controlled dynamical system \eqref{eq:StateEq}
with the control signal acting as an inhomogeneity, the approach of
exactly realizable trajectories leads to Eqs. \eqref{eq:OutputConstraintDAE}
with the desired output $\boldsymbol{y}_{d}$ acting as an inhomogeneity.
As the examples of the last section demonstrate, Eqs. \eqref{eq:OutputConstraintDAE}
are often considerably simpler than the original controlled dynamical
system. In both cases, no simple analytical solution in closed form
can be found for the controlled dynamical system. In general, if the
constraint equation as well as the output relation is linear, the
entire controlled system can be regarded, in some sense and to some
extent, as being linear, even though the original system \eqref{eq:StateEq}
is nonlinear. Three conditions must be met for Eqs. \eqref{eq:OutputConstraintDAE}
to be linear.\newdefinition{asm}{Assumption}
\begin{asm}{\textit{Linearizing assumption.} }\textit{\label{asm:LinearizingAssumption}}
\begin{enumerate}
\item The projectors $\boldsymbol{\mathcal{P}}\left(\boldsymbol{x}\right)$
and $\boldsymbol{\mathcal{Q}}\left(\boldsymbol{x}\right)$ must be
independent of the state $\boldsymbol{x}$,
\begin{align}
\boldsymbol{\mathcal{P}}\left(\boldsymbol{x}\right) & =\boldsymbol{\mathcal{P}}=\boldsymbol{1}-\boldsymbol{\mathcal{Q}}=\text{const},\label{eq:LinearizingAssumption1}
\end{align}

\item the nonlinearity $\boldsymbol{R}\left(\boldsymbol{x}\right)$ must
be linear with respect to the input matrix $\boldsymbol{\mathcal{B}}\left(\boldsymbol{x}\right)$
in the sense that 
\begin{align}
\boldsymbol{\mathcal{Q}}\boldsymbol{R}\left(\boldsymbol{x}\right) & =\boldsymbol{\mathcal{Q}}\boldsymbol{\mathcal{A}}\boldsymbol{x}+\boldsymbol{\mathcal{Q}}\boldsymbol{b}\label{eq:LinearizingAssumption2}
\end{align}
with $n\times n$ matrix $\boldsymbol{\mathcal{A}}$ and $n$-component
vector $\boldsymbol{b}$ independent of the state $\boldsymbol{x}$,
\item the output relation must be linear,
\begin{align}
\boldsymbol{y}\left(t\right) & =\boldsymbol{\mathcal{C}}\boldsymbol{x}\left(t\right).\label{eq:LinearizingAssumption3}
\end{align}

\end{enumerate}
\end{asm}Employing these three assumptions together with the state
separation Eq. \eqref{eq:SplitStateMN} yields
\begin{align}
\boldsymbol{\mathcal{Q}}\boldsymbol{\mathcal{N}}\boldsymbol{\dot{x}}_{d} & =\boldsymbol{\mathcal{Q}}\boldsymbol{\mathcal{A}}\boldsymbol{\mathcal{N}}\boldsymbol{x}_{d}+\boldsymbol{\mathcal{Q}}\boldsymbol{b}+\boldsymbol{\mathcal{Q}}\boldsymbol{\mathcal{A}}\boldsymbol{\mathcal{C}}^{+}\boldsymbol{y}_{d}-\boldsymbol{\mathcal{Q}}\boldsymbol{\mathcal{C}}^{+}\boldsymbol{\dot{y}}_{d}.\label{eq:OutputConstraintEquation}
\end{align}
This is a system of $n-p-r$ algebraic and $r$ differential equations
for the $n-m$ independent components of $\boldsymbol{\mathcal{N}}\boldsymbol{x}_{d}\left(t\right)$.
The number $r$ is given by the rank of the matrix product $\boldsymbol{\mathcal{Q}}\boldsymbol{\mathcal{N}}$,
\begin{align}
r & =\text{rank}\left(\boldsymbol{\mathcal{Q}}\boldsymbol{\mathcal{N}}\right)\leq\min\left(\text{rank}\left(\boldsymbol{\mathcal{Q}}\right),\text{rank}\left(\boldsymbol{\mathcal{N}}\right)\right)\nonumber \\
 & =\min\left(n-p,n-m\right).
\end{align}
Being an inhomogeneous linear differential-algebraic equation, closed
form solutions for Eq. \eqref{eq:OutputConstraintEquation} can be
obtained, see \cite{Campbell1980Singular,Campbell1982Singular,mehrmann2006differentialalgebraic}.
However, a complete discussion is outside the scope of this article
and reserved for later investigations. Equation \eqref{eq:OutputConstraintEquation}
assumes a particularly simple form if $m=p$ and $r=n-p$. This is
the case for the output $\boldsymbol{\mathcal{C}}=\boldsymbol{\mathcal{\hat{P}}}^{T}$,
or, as long as the input matrix is constant, $\boldsymbol{\mathcal{C}}=\boldsymbol{\mathcal{B}}^{T}$.
These outputs imply $\boldsymbol{\mathcal{M}}=\boldsymbol{\mathcal{P}}$
and $\boldsymbol{\mathcal{N}}=\boldsymbol{\mathcal{Q}}$, and Eq.
\eqref{eq:OutputConstraintEquation} becomes a system of $n-p$ differential
equations for the $n-p$ independent components of $\boldsymbol{\mathcal{Q}}\boldsymbol{x}_{d}$,
\begin{align}
\boldsymbol{\mathcal{Q}}\boldsymbol{\dot{x}}_{d} & =\boldsymbol{\mathcal{Q}}\boldsymbol{\mathcal{A}}\boldsymbol{\mathcal{Q}}\boldsymbol{x}_{d}+\boldsymbol{\mathcal{Q}}\boldsymbol{\mathcal{A}}\boldsymbol{\mathcal{P}}\boldsymbol{x}_{d}+\boldsymbol{\mathcal{Q}}\boldsymbol{b}.
\end{align}

\begin{rem}Note that assumption \eqref{eq:LinearizingAssumption1} does neither
imply that the input matrix $\boldsymbol{\mathcal{B}}\left(\boldsymbol{x}\right)$
nor its Moore-Penrose pseudo inverse $\boldsymbol{\mathcal{B}}^{+}\left(\boldsymbol{x}\right)$
is independent of $\boldsymbol{x}$. Furthermore, while the constraint
equation is linear, the control signal $\boldsymbol{u}=\boldsymbol{\mathcal{B}}^{+}\left(\boldsymbol{x}_{d}\right)\left(\boldsymbol{\dot{x}}_{d}-\boldsymbol{R}\left(\boldsymbol{x}_{d}\right)\right)$
may still depend nonlinearly on the desired state $\boldsymbol{x}_{d}$.\end{rem}\begin{rem}Condition
\eqref{eq:LinearizingAssumption2} is very restrictive. It enforces
$n-p$ components to depend only linearly on the state. However, some
important models of nonlinear dynamics satisfy the linearizing assumption.
Among these are all one-dimensional mechanical control systems discussed
in Example \ref{exp:FHN}. Another prominent example is the FitzHugh-Nagumo
model \cite{fitzhugh1955mathematical}, a prototype model for excitable
systems as e.g. the neuron, with a control acting on the activator
variable \cite{loeber2016optimal}.\end{rem}\begin{rem}Systems
satisfying the linearizing assumption are feedback linearizable \cite{khalil2002nonlinear}
without a state transform. This can be seen as follows. Let the feedback-controlled
system be 
\begin{align}
\boldsymbol{\dot{x}}\left(t\right) & =\boldsymbol{R}\left(\boldsymbol{x}\left(t\right)\right)+\boldsymbol{\mathcal{B}}\left(\boldsymbol{x}\left(t\right)\right)\boldsymbol{u}\left(\boldsymbol{x}\left(t\right)\right).
\end{align}
The control signal $\boldsymbol{u}\in\mathbb{R}^{p}$ is transformed
to the new control signal $\boldsymbol{v}\in\mathbb{R}^{p}$ with
the help of the $n\times p$ matrix $\boldsymbol{\mathcal{H}}$ as
\begin{align}
\boldsymbol{u}\left(\boldsymbol{x}\right) & =-\boldsymbol{\mathcal{B}}^{+}\left(\boldsymbol{x}\right)\left(\boldsymbol{R}\left(\boldsymbol{x}\right)-\boldsymbol{\mathcal{H}}\boldsymbol{v}\left(\boldsymbol{x}\right)\right)
\end{align}
to obtain
\begin{align}
\boldsymbol{\dot{x}}\left(t\right) & =\boldsymbol{\mathcal{Q}}\boldsymbol{\mathcal{A}}\boldsymbol{x}\left(t\right)+\boldsymbol{\mathcal{Q}}\boldsymbol{b}+\boldsymbol{\mathcal{P}}\boldsymbol{\mathcal{H}}\boldsymbol{v}\left(\boldsymbol{x}\left(t\right)\right),
\end{align}
where the linearizing assumption has been applied. Let the projector
$\boldsymbol{\mathcal{P}}$ with rank $p$ be rank-decomposed as $\boldsymbol{\mathcal{P}}=\boldsymbol{\mathcal{F}}\boldsymbol{\mathcal{G}}$
with constant $n\times p$ matrix $\boldsymbol{\mathcal{F}}$ and
constant $p\times n$ matrix $\boldsymbol{\mathcal{G}}$. We introduce
the new control signal $\boldsymbol{\tilde{v}}\left(\boldsymbol{x}\right)=\boldsymbol{\mathcal{G}}\boldsymbol{\mathcal{H}}\boldsymbol{v}\left(\boldsymbol{x}\right)$
to get a linear feedback-controlled system,
\begin{align}
\boldsymbol{\dot{x}}\left(t\right) & =\boldsymbol{\mathcal{Q}}\boldsymbol{\mathcal{A}}\boldsymbol{x}\left(t\right)+\boldsymbol{\mathcal{Q}}\boldsymbol{b}+\boldsymbol{\mathcal{F}}\boldsymbol{\tilde{v}}\left(\boldsymbol{x}\left(t\right)\right).
\end{align}
\end{rem}

\section{\label{sec:Controllability}Controllability}

A system is called controllable, or full state controllable, if it
is possible to achieve a transfer from an initial state $\boldsymbol{x}\left(t_{0}\right)=\boldsymbol{x}_{0}$
at time $t=t_{0}$ to a final state $\boldsymbol{x}\left(t_{1}\right)=\boldsymbol{x}_{1}$
at the terminal time $t=t_{1}$ \cite{kalman1959general}. Along which
trajectory the transfer is achieved is irrelevant. While for LTI systems
conditions for controllability are easily expressed in terms of a
Kalman rank condition, these conditions are more difficult for nonlinear
control systems \cite{khalil2002nonlinear,isidori1995nonlinear}.
Here, we derive a similar rank condition within the framework of exactly
realizable trajectories. This rank condition also applies to systems
satisfying the linearizing assumption \ref{asm:LinearizingAssumption}.\\
We consider the controlled state equation \eqref{eq:StateEq} together
with the linearizing assumption Eq. \eqref{eq:LinearizingAssumption2}.
This implies a linear constraint equation \eqref{eq:ConstraintEquationForRealizableTrajectories},
\begin{align}
\boldsymbol{\mathcal{Q}}\boldsymbol{\dot{x}}_{d} & =\boldsymbol{\mathcal{Q}}\boldsymbol{\mathcal{A}}\boldsymbol{\mathcal{Q}}\boldsymbol{x}_{d}+\boldsymbol{\mathcal{Q}}\boldsymbol{\mathcal{A}}\boldsymbol{\mathcal{P}}\boldsymbol{x}_{d}+\boldsymbol{\mathcal{Q}}\boldsymbol{b}.\label{eq:LinearConstraintEquation}
\end{align}
Because the trajectory in between $\boldsymbol{x}_{0}$ and $\boldsymbol{x}_{1}$
is irrelevant, we may assume that the state components $\boldsymbol{\mathcal{P}}\boldsymbol{x}_{d}\left(t\right)$
are prescribed while the components $\boldsymbol{\mathcal{Q}}\boldsymbol{x}_{d}\left(t\right)$
are governed by Eq. \eqref{eq:LinearConstraintEquation}. Equation
\eqref{eq:LinearConstraintEquation} is an inhomogeneous linear dynamical
system for $\boldsymbol{\mathcal{Q}}\boldsymbol{x}_{d}\left(t\right)$.
Its solution with initial condition $\boldsymbol{\mathcal{Q}}\boldsymbol{x}_{d}\left(t_{0}\right)=\boldsymbol{\mathcal{Q}}\boldsymbol{x}_{0}$
yields 
\begin{align}
\boldsymbol{\mathcal{Q}}\boldsymbol{x}_{d}\left(t\right) & =\int_{t_{0}}^{t}\text{d}\tau\exp\left(\boldsymbol{\mathcal{Q}}\boldsymbol{\mathcal{A}}\boldsymbol{\mathcal{Q}}\left(t-\tau\right)\right)\boldsymbol{\mathcal{Q}}\left(\boldsymbol{\mathcal{A}}\boldsymbol{\mathcal{P}}\boldsymbol{x}_{d}\left(\tau\right)+\boldsymbol{b}\right)\nonumber \\
 & +\exp\left(\boldsymbol{\mathcal{Q}}\boldsymbol{\mathcal{A}}\boldsymbol{\mathcal{Q}}\left(t-t_{0}\right)\right)\boldsymbol{\mathcal{Q}}\boldsymbol{x}_{0}.\label{eq:Qx_dSolution}
\end{align}

\begin{thm}{Controllability.\newline}\label{thm:Controllability}A nonlinear system Eq. \eqref{eq:StateEq}
which satisfies the linearizing assumption \ref{asm:LinearizingAssumption}
is controllable if the $n\times n^{2}$ controllability matrix
\begin{align}
\boldsymbol{\mathcal{K}} & =\left(\boldsymbol{\mathcal{Q}}\boldsymbol{\mathcal{A}}\boldsymbol{\mathcal{P}}|\boldsymbol{\mathcal{Q}}\boldsymbol{\mathcal{A}}\boldsymbol{\mathcal{Q}}\boldsymbol{\mathcal{A}}\boldsymbol{\mathcal{P}}|\cdots|\left(\boldsymbol{\mathcal{Q}}\boldsymbol{\mathcal{A}}\boldsymbol{\mathcal{Q}}\right)^{n-1}\boldsymbol{\mathcal{Q}}\boldsymbol{\mathcal{A}}\boldsymbol{\mathcal{P}}\right).\label{eq:ControllabilityMatrixForRealizableTrajectories}
\end{align}
satisfies the rank condition 
\begin{align}
\text{\textup{rank}}\left(\boldsymbol{\mathcal{K}}\right) & =n-p.\label{eq:RankConditionRealizableTrajectories}
\end{align}
\end{thm}
\newproof{pot2}{Proof of Theorem \ref{thm:Controllability}}
\begin{pot2}\label{pot:pot2}Achieving a transfer from an initial to a finite
state means the desired trajectory $\boldsymbol{x}_{d}\left(t\right)$
has to satisfy $\boldsymbol{x}_{d}\left(t_{0}\right)=\boldsymbol{x}_{0}$
and $\boldsymbol{x}_{d}\left(t_{1}\right)=\boldsymbol{x}_{1}$. Consequently,
$\boldsymbol{\mathcal{P}}\boldsymbol{x}_{d}$ and $\boldsymbol{\mathcal{Q}}\boldsymbol{x}_{d}$
have to satisfy the initial and terminal conditions
\begin{align}
\boldsymbol{\mathcal{P}}\boldsymbol{x}_{d}\left(t_{0}\right) & =\boldsymbol{\mathcal{P}}\boldsymbol{x}_{0}, & \boldsymbol{\mathcal{P}}\boldsymbol{x}_{d}\left(t_{1}\right) & =\boldsymbol{\mathcal{P}}\boldsymbol{x}_{1},\label{eq:Px_dInitialAndTerminalConstraints}\\
\boldsymbol{\mathcal{Q}}\boldsymbol{x}_{d}\left(t_{0}\right) & =\boldsymbol{\mathcal{Q}}\boldsymbol{x}_{0}, & \boldsymbol{\mathcal{Q}}\boldsymbol{x}_{d}\left(t_{1}\right) & =\boldsymbol{\mathcal{Q}}\boldsymbol{x}_{1}.\label{eq:Qx_dInitialAndTerminalConstraints}
\end{align}
The part $\boldsymbol{\mathcal{P}}\boldsymbol{x}_{d}$ is prescribed
by the experimenter such that it satisfies the initial and terminal
conditions Eq. \eqref{eq:Px_dInitialAndTerminalConstraints}. Consequently,
all initial and terminal conditions except $\boldsymbol{\mathcal{Q}}\boldsymbol{x}_{d}\left(t_{1}\right)=\boldsymbol{\mathcal{Q}}\boldsymbol{x}_{1}$
are satisfied. Enforcing this remaining condition onto the solution
Eq. \eqref{eq:Qx_dSolution} of the constraint equation yields
\begin{align}
\boldsymbol{\mathcal{Q}}\boldsymbol{x}_{1} & =\boldsymbol{\mathcal{Q}}\boldsymbol{x}_{d}\left(t_{1}\right)=\exp\left(\boldsymbol{\mathcal{Q}}\boldsymbol{\mathcal{A}}\boldsymbol{\mathcal{Q}}\left(t_{1}-t_{0}\right)\right)\boldsymbol{\mathcal{Q}}\boldsymbol{x}_{0}\label{eq:Px_dConstraint}\\
 & +\int_{t_{0}}^{t_{1}}\text{d}\tau\exp\left(\boldsymbol{\mathcal{Q}}\boldsymbol{\mathcal{A}}\boldsymbol{\mathcal{Q}}\left(t_{1}-\tau\right)\right)\boldsymbol{\mathcal{Q}}\left(\boldsymbol{\mathcal{A}}\boldsymbol{\mathcal{P}}\boldsymbol{x}_{d}\left(\tau\right)+\boldsymbol{b}\right).\nonumber 
\end{align}
This can actually be viewed as a condition for the part $\boldsymbol{\mathcal{P}}\boldsymbol{x}_{d}$.
The transfer from $\boldsymbol{x}_{0}$ to $\boldsymbol{x}_{1}$ is
achieved as long as $\boldsymbol{\mathcal{P}}\boldsymbol{x}_{d}$
satisfies Eq. \eqref{eq:Px_dConstraint}.

Similarly to the proof of the Kalman rank condition in \cite{chen1995linear},
conditions on the state matrix $\boldsymbol{\mathcal{A}}$ and the
projectors $\boldsymbol{\mathcal{P}}$ and $\boldsymbol{\mathcal{Q}}$
can be given such that Eq. \eqref{eq:Px_dConstraint} is satisfied.
Due to the Cayley-Hamilton theorem \cite{axler1997linear}, any power
of matrices with $i\geq n$ can be expanded in terms of lower order
matrix powers as 
\begin{align}
\left(\boldsymbol{\mathcal{Q}}\boldsymbol{\mathcal{A}}\boldsymbol{\mathcal{Q}}\right)^{i} & =\sum_{k=0}^{n-1}d_{ik}\left(\boldsymbol{\mathcal{Q}}\boldsymbol{\mathcal{A}}\boldsymbol{\mathcal{Q}}\right)^{k}.\label{eq:CayleyHamiltonForQAQ}
\end{align}
The term involving $\boldsymbol{\mathcal{P}}\boldsymbol{x}_{d}\left(\tau\right)$
in Eq. \eqref{eq:Px_dConstraint} can be simplified,
\begin{align}
 & \int_{t_{0}}^{t_{1}}\text{d}\tau\exp\left(\boldsymbol{\mathcal{Q}}\boldsymbol{\mathcal{A}}\boldsymbol{\mathcal{Q}}\left(t_{0}-\tau\right)\right)\boldsymbol{\mathcal{Q}}\boldsymbol{\mathcal{A}}\boldsymbol{\mathcal{P}}\boldsymbol{x}_{d}\left(\tau\right)\nonumber \\
= & \sum_{k=0}^{n-1}\left(\boldsymbol{\mathcal{Q}}\boldsymbol{\mathcal{A}}\boldsymbol{\mathcal{Q}}\right)^{k}\boldsymbol{\mathcal{Q}}\boldsymbol{\mathcal{A}}\boldsymbol{\mathcal{P}}\nonumber \\
 & \times\int_{t_{0}}^{t_{1}}\text{d}\tau\left(\frac{\left(t_{0}-\tau\right)^{k}}{k!}+\sum_{i=n}^{\infty}d_{ik}\frac{\left(t_{0}-\tau\right)^{i}}{i!}\right)\boldsymbol{\mathcal{P}}\boldsymbol{x}_{d}\left(\tau\right),
\end{align}
such that Eq. \eqref{eq:Px_dConstraint} becomes a truncated sum
\begin{align}
 & \exp\left(-\boldsymbol{\mathcal{Q}}\boldsymbol{\mathcal{A}}\boldsymbol{\mathcal{Q}}\left(t_{1}-t_{0}\right)\right)\boldsymbol{\mathcal{Q}}\boldsymbol{x}_{1}\nonumber \\
 & -\boldsymbol{\mathcal{Q}}\boldsymbol{x}_{0}-\int_{t_{0}}^{t_{1}}\text{d}\tau\exp\left(\boldsymbol{\mathcal{Q}}\boldsymbol{\mathcal{A}}\boldsymbol{\mathcal{Q}}\left(t_{0}-\tau\right)\right)\boldsymbol{\mathcal{Q}}\boldsymbol{b}\nonumber \\
= & \sum_{k=0}^{n-1}\left(\boldsymbol{\mathcal{Q}}\boldsymbol{\mathcal{A}}\boldsymbol{\mathcal{Q}}\right)^{k}\boldsymbol{\mathcal{Q}}\boldsymbol{\mathcal{A}}\boldsymbol{\mathcal{P}}\boldsymbol{\alpha}_{k}\left(t_{1},t_{0}\right).\label{eq:Px_rConstraint_2}
\end{align}
We defined the $n\times1$ vectors $\boldsymbol{\alpha}_{k}$ for
$k\in\left\{ 0,\dots,n-1\right\} $ as
\begin{align}
 & \boldsymbol{\alpha}_{k}\left(t_{1},t_{0}\right)\nonumber \\
 & =\int_{t_{0}}^{t_{1}}\text{d}\tau\left(\frac{\left(t_{0}-\tau\right)^{k}}{k!}+\sum_{i=n}^{\infty}d_{ik}\frac{\left(t_{0}-\tau\right)^{i}}{i!}\right)\boldsymbol{\mathcal{P}}\boldsymbol{x}_{d}\left(\tau\right).
\end{align}
The right hand side of Eq. \eqref{eq:Px_rConstraint_2} can be written
with the help of the $n^{2}\times1$ vector 
\begin{align}
\boldsymbol{\alpha}\left(t_{1},t_{0}\right) & =\left(\boldsymbol{\alpha}_{0}\left(t_{1},t_{0}\right),\dots,\boldsymbol{\alpha}_{n-1}\left(t_{1},t_{0}\right)\right)^{T}
\end{align}
as 
\begin{align}
\exp\left(-\boldsymbol{\mathcal{Q}}\boldsymbol{\mathcal{A}}\boldsymbol{\mathcal{Q}}\left(t_{1}-t_{0}\right)\right)\boldsymbol{\mathcal{Q}}\boldsymbol{x}_{1}-\boldsymbol{\mathcal{Q}}\boldsymbol{x}_{0}\nonumber \\
-\int_{t_{0}}^{t_{1}}\text{d}\tau\exp\left(\boldsymbol{\mathcal{Q}}\boldsymbol{\mathcal{A}}\boldsymbol{\mathcal{Q}}\left(t_{0}-\tau\right)\right)\boldsymbol{\mathcal{Q}}\boldsymbol{b} & =\boldsymbol{\mathcal{K}}\boldsymbol{\alpha}\left(t_{1},t_{0}\right).\label{eq:Px_rConstraint_3}
\end{align}
We defined the $n\times n^{2}$ controllability matrix \textit{$\boldsymbol{\mathcal{K}}$}
as given by Eq. \eqref{eq:ControllabilityMatrixForRealizableTrajectories}.
The left hand side of Eq. \eqref{eq:Px_rConstraint_3} can be any
point in $\boldsymbol{\mathcal{Q}}\mathbb{R}^{n}=\mathbb{R}^{n-p}$.
The mapping from $\boldsymbol{\mathcal{Q}}\boldsymbol{x}_{1}$ to
$\boldsymbol{\alpha}$ is surjective, i.e., every element on the left
hand side has a corresponding element on the right hand side, if $\boldsymbol{\mathcal{K}}$
has full rank $n-p$.\qed\end{pot2}
\begin{rem}Using the complementary projectors $\boldsymbol{\mathcal{P}}$ and
$\boldsymbol{\mathcal{Q}}$, the state matrix $\boldsymbol{\mathcal{A}}$
can be split up in four parts as 
\begin{align}
\boldsymbol{\mathcal{A}} & =\boldsymbol{\mathcal{P}}\boldsymbol{\mathcal{A}}\boldsymbol{\mathcal{P}}+\boldsymbol{\mathcal{P}}\boldsymbol{\mathcal{A}}\boldsymbol{\mathcal{Q}}+\boldsymbol{\mathcal{Q}}\boldsymbol{\mathcal{A}}\boldsymbol{\mathcal{P}}+\boldsymbol{\mathcal{Q}}\boldsymbol{\mathcal{A}}\boldsymbol{\mathcal{Q}}.
\end{align}
Note that the controllability matrix $\boldsymbol{\mathcal{\tilde{K}}}$,
Eq. \eqref{eq:ControllabilityMatrixForRealizableTrajectories}, does
only depend on the parts $\boldsymbol{\mathcal{Q}}\boldsymbol{\mathcal{A}}\boldsymbol{\mathcal{P}}$
and $\boldsymbol{\mathcal{Q}}\boldsymbol{\mathcal{A}}\boldsymbol{\mathcal{Q}}$,
but not on $\boldsymbol{\mathcal{P}}\boldsymbol{\mathcal{A}}\boldsymbol{\mathcal{P}}$
and $\boldsymbol{\mathcal{P}}\boldsymbol{\mathcal{A}}\boldsymbol{\mathcal{Q}}$.
Consequently, only knowledge of the parts $\boldsymbol{\mathcal{Q}}\boldsymbol{\mathcal{A}}\boldsymbol{\mathcal{P}}$
and $\boldsymbol{\mathcal{Q}}\boldsymbol{\mathcal{A}}\boldsymbol{\mathcal{Q}}$
is required to decide if a system is controllable. This might be advantageous
for applications with incomplete knowledge about the underlying dynamics.\end{rem}

\section{\label{sec:OutputControllability}Output controllability}

\begin{thm}{Output controllability.\newline}\label{thm:Theorem5}The system \eqref{eq:StateEq} satisfying the
linearizing assumption \ref{asm:LinearizingAssumption} is output
controllable if the output controllability matrix
\begin{align}
\boldsymbol{\mathcal{K}}_{\boldsymbol{\mathcal{C}}} & =\left(\boldsymbol{\mathcal{C}}\boldsymbol{\mathcal{P}}|\boldsymbol{\mathcal{C}}\boldsymbol{\mathcal{Q}}\boldsymbol{\mathcal{A}}\boldsymbol{\mathcal{P}}|\cdots|\boldsymbol{\mathcal{C}}\left(\boldsymbol{\mathcal{Q}}\boldsymbol{\mathcal{A}}\boldsymbol{\mathcal{Q}}\right)^{n-1}\boldsymbol{\mathcal{Q}}\boldsymbol{\mathcal{A}}\boldsymbol{\mathcal{P}}\right)\label{eq:OutputControllabilityMatrix}
\end{align}
satisfies the rank condition
\begin{align}
\text{\textup{rank}}\left(\boldsymbol{\mathcal{K}}_{\boldsymbol{\mathcal{C}}}\right) & =m.
\end{align}
\end{thm}
\newproof{pot5}{Proof of Theorem \ref{thm:Theorem5}}
\begin{pot5}Using Eq. \eqref{eq:Qx_dSolution}, the solution for the linear output
reads
\begin{align}
\boldsymbol{y}_{d}\left(t\right) & =\boldsymbol{\mathcal{C}}\boldsymbol{\mathcal{P}}\boldsymbol{x}_{d}\left(t\right)+\boldsymbol{\mathcal{C}}\exp\left(\boldsymbol{\mathcal{Q}}\boldsymbol{\mathcal{A}}\boldsymbol{\mathcal{Q}}\left(t-t_{0}\right)\right)\boldsymbol{\mathcal{Q}}\boldsymbol{x}_{0}\nonumber \\
 & +\boldsymbol{\mathcal{C}}\int_{t_{0}}^{t}\text{d}\tau\exp\left(\boldsymbol{\mathcal{Q}}\boldsymbol{\mathcal{A}}\boldsymbol{\mathcal{Q}}\left(t-\tau\right)\right)\boldsymbol{\mathcal{Q}}\boldsymbol{\mathcal{A}}\boldsymbol{\mathcal{P}}\boldsymbol{x}_{d}\left(\tau\right).\label{eq:Eq88}
\end{align}
Exploiting the Cayley-Hamilton theorem, Eq. \eqref{eq:Eq88} becomes
\begin{align}
\boldsymbol{y}_{d}\left(t\right)-\boldsymbol{\mathcal{C}}\exp\left(\boldsymbol{\mathcal{Q}}\boldsymbol{\mathcal{A}}\boldsymbol{\mathcal{Q}}\left(t-t_{0}\right)\right)\boldsymbol{\mathcal{Q}}\boldsymbol{x}_{0}\nonumber \\
=\boldsymbol{\mathcal{C}}\boldsymbol{\mathcal{P}}\boldsymbol{x}_{d}\left(t\right)+\boldsymbol{\mathcal{C}}\sum_{k=0}^{n-1}\left(\boldsymbol{\mathcal{Q}}\boldsymbol{\mathcal{A}}\boldsymbol{\mathcal{Q}}\right)^{k}\boldsymbol{\mathcal{Q}}\boldsymbol{\mathcal{A}}\boldsymbol{\mathcal{P}}\boldsymbol{\beta}_{k}\left(t,t_{0}\right).
\end{align}
The remainder of the proof proceeds analogously to the proof of Theorem
\ref{thm:Controllability} for full state controllability.\qed\end{pot5}
\begin{rem}For $\boldsymbol{\mathcal{C}}=\boldsymbol{1}$ and $m=n$, output
controllability should reduce to full state controllability. Indeed,
if the controllability matrix $\boldsymbol{\mathcal{K}}$, Eq. \eqref{eq:ControllabilityMatrixForRealizableTrajectories},
satisfies the rank condition $\text{rank}\left(\boldsymbol{\mathcal{K}}\right)=n-p$,
then the matrix $\boldsymbol{\mathcal{\tilde{K}}}=\left(\boldsymbol{\mathcal{P}}|\boldsymbol{\mathcal{K}}\right)$
satisfies $\text{rank}\left(\boldsymbol{\mathcal{\tilde{K}}}\right)=n$,
and $\boldsymbol{\mathcal{K}}_{\boldsymbol{\mathcal{C}}}$ reduces
to $\boldsymbol{\mathcal{\tilde{K}}}$ for $\boldsymbol{\mathcal{C}}=\boldsymbol{1}$
and $m=n$.\end{rem}\begin{rem}Similar as for controllability, we
can use the Moore-Penrose projectors $\boldsymbol{\mathcal{M}}$ and
$\boldsymbol{\mathcal{N}}$ constructed from the output matrix $\boldsymbol{\mathcal{C}}$
in Definition \ref{def:Definition2} to express observability in terms
of a rank condition for an observability matrix \cite{loeber2016optimal}.\end{rem}

\section{\label{sec:2Conclusions}Conclusions and outlook}

Exactly realizable desired trajectories are the subset of desired
trajectories $\boldsymbol{x}_{d}\left(t\right)$ for which a control
exists such that the state over time $\boldsymbol{x}\left(t\right)$
follows the desired trajectory exactly, $\boldsymbol{x}\left(t\right)=\boldsymbol{x}_{d}\left(t\right)$.
By means of the Moore-Penrose projectors defined in Eqs. \eqref{eq:TimeDependentPQ},
we propose a separation of the state equation \eqref{eq:StateEq}
in two parts. The first part, called the constraint equation \eqref{eq:ConstraintEquationForRealizableTrajectories},
is independent of the control signal $\boldsymbol{u}$. The second
part Eq. \eqref{eq:ControlSignal} establishes a one-to-one relationship
between the $p$-dimensional control signal $\boldsymbol{u}\left(t\right)$
and $p$ out of $n$ components of the desired trajectory $\boldsymbol{x}_{d}\left(t\right)$.
The constraint equation fixes those $n-p$ components of the desired
trajectory $\boldsymbol{x}_{d}\left(t\right)$ for which no one-to-one
relationship with the control signal exists. A desired trajectory
is exactly realizable if and only if it satisfies the constraint equation.\\
We can distinguish 3 classes of desired trajectories $\boldsymbol{x}_{d}$:
\begin{enumerate}[(A)]
\item desired trajectories $\boldsymbol{x}_{d}$ which are solutions to
the uncontrolled system,
\item desired trajectories $\boldsymbol{x}_{d}$ which are exactly realizable,
\item arbitrary desired trajectories $\boldsymbol{x}_{d}$.\end{enumerate}
Desired trajectories of class (A) satisfy the uncontrolled state equation
\begin{align}
\boldsymbol{\dot{x}}_{d}\left(t\right) & =\boldsymbol{R}\left(\boldsymbol{x}_{d}\left(t\right)\right).\label{eq:Eq2249}
\end{align}
This constitutes the most specific class of desired trajectories.
Because of Eq. \eqref{eq:Eq2249}, the constraint equation \eqref{eq:ConstraintEquationForRealizableTrajectories}
is trivially satisfied and the control signal given by Eq. \eqref{eq:ControlSignal}
vanishes, $\boldsymbol{u}\left(t\right)=\boldsymbol{0}$. Class (A)
encompasses several important control tasks, as e.g. the stabilization
of unstable stationary states and periodic orbits \cite{ott1990controlling,sontag2011stability}.
Only for desired trajectories of class (A) it is possible to find
non-invasive controls. Non-invasive control signals vanish upon achieving
the control target. The open-loop control approach developed here
cannot be employed to trajectories of class (A). Instead, the stabilization
of unstable solutions to uncontrolled systems requires feedback control.\\
Desired trajectories of class (B) satisfy the constraint equation
\eqref{eq:ConstraintEquationForRealizableTrajectories} and yield
a non-vanishing control signal $\boldsymbol{u}\left(t\right)\neq\boldsymbol{0}$.
The approach developed here applies to this class. Several other techniques
developed in mathematical control theory as e.g. feedback linearization
and differential flatness, also work with this class of desired trajectories
\cite{khalil2002nonlinear,sira2004differentially}. Class (B) contains
the desired trajectories from class (A) as a special case. For desired
trajectories of class (A) and class (B), the solution of the controlled
state is simply given by $\boldsymbol{x}\left(t\right)=\boldsymbol{x}_{d}\left(t\right)$.\\
Finally, class (C) is the most general class of desired trajectories
and contains class (A) and (B) as special cases. In general, these
desired trajectories do not satisfy the constraint equation,
\begin{align}
\boldsymbol{0} & \neq\boldsymbol{\mathcal{Q}}\left(\boldsymbol{x}_{d}\left(t\right)\right)\left(\boldsymbol{\dot{x}}_{d}\left(t\right)-\boldsymbol{R}\left(\boldsymbol{x}_{d}\left(t\right)\right)\right),
\end{align}
such that the approach developed here cannot be applied to all desired
trajectories of class (C). A general expression for the control signal
in terms of the desired trajectory $\boldsymbol{x}_{d}\left(t\right)$
is not available. Furthermore, the solution for the controlled state
trajectory $\boldsymbol{x}\left(t\right)$ is usually not simply given
by $\boldsymbol{x}_{d}\left(t\right)$, $\boldsymbol{x}\left(t\right)\neq\boldsymbol{x}_{d}\left(t\right)$.
A solution to control problems defined by class (C) does not only
consist in finding an expression for the control signal, but also
involves finding a solution for the controlled state $\boldsymbol{x}\left(t\right)$.
One possible method to solve such control problems is optimal trajectory
tracking. This technique is concerned with minimizing the distance
between $\boldsymbol{x}\left(t\right)$ and $\boldsymbol{x}_{d}\left(t\right)$
in function space as measured by the functional 
\begin{align}
\mathcal{J} & =\frac{1}{2}\intop_{t_{0}}^{t_{1}}\text{d}t\left(\boldsymbol{x}\left(t\right)-\boldsymbol{x}_{d}\left(t\right)\right)^{2}+\dfrac{\epsilon^{2}}{2}\intop_{t_{0}}^{t_{1}}\text{d}t\left(\boldsymbol{u}\left(t\right)\right)^{2}.
\end{align}
The functional $\mathcal{J}$ is to be minimized subject to the constraint
that $\boldsymbol{x}\left(t\right)$ is given as the solution to the
controlled dynamical system \eqref{eq:StateEq}. The regularization
term with small coefficient $0<\epsilon\ll1$ ensures the existence
of solutions for $\boldsymbol{x}\left(t\right)$ and $\boldsymbol{u}\left(t\right)$
within appropriate function spaces \cite{bryson1969applied}. For
$\epsilon\rightarrow0$, the state as well as the control may diverge
\cite{loeber2016optimal} and the optimization procedure becomes a
so-called singular optimal control problem \cite{bryson1969applied}.
However, because exactly realizable desired trajectories satisfy $\boldsymbol{x}\left(t\right)=\boldsymbol{x}_{d}\left(t\right)$,
they can be viewed as bounded solutions to unregularized ($\epsilon=0$)
optimal trajectory tracking problems.

The linearizing assumption defines a class of nonlinear control systems
which essentially behave like linear control system. Systems satisfying
the linearizing assumption allow exact analytical solutions in closed
form even if no analytical solutions for the uncontrolled system exist.
The linearizing assumption uncovers a hidden linear structure underlying
nonlinear open-loop control systems. Similarly, feedback linearization
defines a huge class of nonlinear control systems possessing an underlying
linear structure \cite{khalil2002nonlinear}. The class of feedback
linearizable systems contains the systems satisfying the linearizing
assumption as a trivial case. However, the linearizing assumption
defined here goes much further than feedback linearization. We were
able to apply the relatively simple notion of controllability in terms
of a rank condition to systems satisfying the linearizing assumption.
This is a direct extension of the properties of linear control systems
to a class of nonlinear control systems. Furthermore, we may combine
the linearizing assumption with the viewpoint that exactly realizable
trajectories solve an unregularized optimal control problem. This
reveals the possibility of linear structures underlying nonlinear
optimal trajectory tracking in the limit of vanishing regularization
parameter $\epsilon\rightarrow0$ \cite{loeber2016optimal}.

Finally, we mention a possible extension of the ideas expounded here
to spatio-temporal systems. While generalizing the notion of controllability
to spatiotemporal systems encounters difficulties due to an infinite-dimensional
state space, generalizing the notion of an exactly realizable trajectory
is straightforward \cite{loeber2016optimal}. We applied these ideas
in a slightly different form to control the position, orientation,
and shape of wave patterns in reaction-diffusion systems in \cite{lober2014controlling,lober2014shaping,lober2014stability,lober2014control}.

\bibliographystyle{elsarticle-num}
\bibliography{ExactlyRealizableTrajectories}

\end{document}